\documentclass[12pt,reqno,a4]{amsart}

\usepackage{a4}
\usepackage{graphicx}
\usepackage{amsmath,amssymb,natbib,color,bbm,ifthen}
\usepackage[latin1]{inputenc}

\def\BIBDIR{../bib}

\newtheorem{theorem}{Theorem}
\newtheorem{lemma}[theorem]{Lemma}
\newtheorem{definition}{Definition}
\newtheorem{proposition}[theorem]{Proposition}
\newtheorem{corollary}[theorem]{Corollary}
\theoremstyle{remark}
\newtheorem{remark}{Remark}

  \advance\oddsidemargin -1.0cm
  \advance\evensidemargin -1.0cm

\renewcommand{\hat}{\widehat}

\newenvironment{enum_W}
  {%
  \setlength{\leftmargini}{4em}\begin{enumerate}}
  {\end{enumerate}}

  {\end{enumerate}}

\newenvironment{enum_a}
  {%
  \setlength{\leftmargini}{2em}\begin{enumerate}}%
  {\end{enumerate}}

\def\rme{\mathrm{e}}

\def\dj{u}
\def\dk{v}
\def\rme{\mathrm{e}}
\def\rmi{\mathrm{i}}
\def\ltwo{\mathrm{L}^2}
\def\allWA{\ref{item:Wreg}-\ref{item:MIM}}
\def\trace{\mathrm{Tr}}

\def\eqdist{\stackrel{\mathrm{d}}{=}}

\def\dwt{W}
\def\bdwt{\mathbf{\dwt}}
\def\indexset{\mathcal{I}}

\def\L{\mathrm{T}}

\def\be{\mathbf{e}}

\def\diffop{\mathbf{\Delta}}

\def\1{\mathbbm{1}}

\def\bg{\mathbf{g}}
\def\bc{\mathbf{c}}
\def\bx{\mathbf{x}}
\def\bY{\mathbf{Y}}
\def\Nset{\mathbb{N}} 
\def\Zset{\mathbb{Z}}

\def\Rset{\mathbb{R}} 
\def\Cset{\mathbb{C}} 
\def\PE{\mathbb{E}} 
\def\PVar{\mathrm{Var}}
\def\PCov{\mathrm{Cov}}

\def\calN{\mathcal{N}}

\def\calH{\mathcal{H}}

\def\ie{\textit{i.e.} }
\def\eg{\textit{e.g.} }

\def\prob{\mathbb{P}}

\newcommand{\eqdef}{\ensuremath{\stackrel{\mathrm{def}}{=}}}
\newcommand{\eqsp}{\;}

\newcommand{\AVvar}[3][]
{
\ifthenelse{\equal{#1}{}}{\mathbf{V}_{#3}(#2)}{\mathbf{V}_{#3}(#2,#1)}}
\newcommand{\AVvarJoint}[3][]
{
\ifthenelse{\equal{#1}{}}{\mathbf{W}_{#3}(#2)}{\mathbf{W}_{#3}(#2,#1)}}

\newcommand{\AsympVarWWE}[2][]
{ \ifthenelse{\equal{#1}{}}{\mathrm{V}(#2)}{\mathrm{V}(#2,#1)}}

\newcommand{\AVvarInv}[3][]
{\ifthenelse{\equal{#1}{}}{\mathbf{V}^{-1}_{#3}(#2)}{\mathbf{V}^{-1}_{#3}(#2,#1)}}
\newcommand{\sigmaasymp}[2][]
{
\ifthenelse{\equal{#1}{}}{\sigma(#2)}{\sigma(#2,#1)}}

\def\dmin{{d_{\min}}}
\def\dmax{{d_{\max}}}
\def\vjsymb{\sigma}
\newcommand{\vj}[3][]{%
\ifthenelse{\equal{#1}{}}{\vjsymb^2_{#2}(#3)}{\vjsymb^2_{#2}(#3,#1)}}
\newcommand{\vjj}[3][]{%
\ifthenelse{\equal{#1}{}}{\vjsymb^4_{#2}(#3)}{\vjsymb^4_{#2}(#3,#1)}}

\newcommand{\hvj}[3][]{%
\ifthenelse{\equal{#1}{}}{\hat{\vjsymb}^2_{#2}}{\hat{\vjsymb}^2_{#2}(#1)}}
\newcommand{\Kvar}[2][]{
\ifthenelse{\equal{#1}{}}{\mathrm{K}(#2)}{\mathrm{K}(#2,#1)}}

\def\densletter{\mathbf{D}}
\newcommand{\bdens}[4][]{%
\ifthenelse{\equal{#1}{}}{\densletter_{#2}({#3};#4)}{\densletter_{#2}(#3;#4,#1)}}

\def\densasympletter{D}
\newcommand{\densasymp}[3][]{%
\ifthenelse{\equal{#1}{}}{\mathrm{\densasympletter}(#2;#3)}{\mathrm{\densasympletter}(#2;#3,#1)}}
\newcommand{\bdensasymp}[4][]{%
\ifthenelse{\equal{#1}{}}{\mathbf{\densasympletter}_{\infty,#2}({#3};#4)}{\mathbf{\densasympletter}_{\infty,#2}(#3;#4,#1)}}

\def\intbdensletter{I}
\newcommand{\intbdens}[3][]{
\ifthenelse{\equal{#1}{}}{\mathrm{\intbdensletter}_{#2}(#3)}{\mathrm{\intbdensletter}_{#2}(#3,#1)}}

\def\regressweights{\mathbf{w}}

\title[Spectral Density of the Wavelet Coefficients]{On the spectral density of the wavelet coefficients of long memory time series with application
to the log-regression estimation of the memory parameter}

\author{E. Moulines}
\author{F. Roueff}

\address{GET/Télécom Paris, CNRS LTCI, 46, rue Barrault, 75634 Paris Cédex 13, France.}
\email{moulines@tsi.enst.fr}
\email{roueff@tsi.enst.fr}
\author{M.S. Taqqu}
\address{Department of Mathematics and Statistics, Boston University Boston, MA 02215, USA.}
\email{murad@math.bu.edu}



\thanks{Murad~S.~Taqqu would like to thank l'\'Ecole Normale Sup\'erieure des T\'elecom\-munications in Paris
for their hospitality.  This research was partially supported by the NSF Grant DMS--0505747 at Boston University.}

\date{June 6, 2006}

\begin{document}

\maketitle

\begin{center}
{\it CNRS LTCI and Boston University}\\
\end{center}
\renewcommand{\thefootnote}{}
\footnote{\textit{Corresponding author}: F. Roueff,  GET/Télécom Paris, CNRS LTCI, 46, rue Barrault, 75634 Paris Cédex 13, France.}
\renewcommand{\thefootnote}{\arabic{footnote}}

\begin{abstract}
In recent years, methods to estimate the memory parameter using wavelet analysis have gained popularity
in many areas of science.
Despite its widespread use, a rigorous semi-parametric asymptotic theory, comparable
to the one developed for Fourier methods, is still missing. In this contribution, we adapt to the wavelet setting the
classical semi-parametric framework introduced by Robinson and his co-authors for estimating the memory parameter
of a (possibly) non-stationary process. 
Our results apply to a class of wavelets with bounded supports, which include but are not limited to Daubechies wavelets. We
derive an explicit expression of the spectral density of the wavelet coefficients and show that it can be approximated, at
large scales, by the spectral density of the continuous-time wavelet coefficients of fractional Brownian motion. We derive an
explicit bound for the difference between the spectral densities.
As an application, we obtain minimax upper bounds for the log-scale regression estimator of the memory parameter
for a Gaussian process and we derive an explicit expression of its asymptotic variance.
\end{abstract}

\newpage

\section{Introduction}
Let  $X \eqdef \{X_k \}_{k\in\Zset}$ be a real-valued process, not necessarily stationary and let $\diffop^K X$ denote its
$K$-th order difference. The first order difference is $[\diffop X]_k \eqdef X_k - X_{k-1}$ and $\diffop^K$ is defined recursively.
The process $X$ is said to have \emph{memory parameter} $d$, $d \in \Rset$ (in short, is an M($d$) process) if for any integer $K>d-1/2$,
the $K$-th order difference process $\diffop^{K}X$ is weakly stationary with spectral density function
\begin{equation}\label{eq:spdelta}
f_{\diffop^K X}(\lambda) \eqdef |1-e^{-i\lambda}|^{2(K-d)}\,f^\ast(\lambda)\quad  \lambda\in(-\pi,\pi),
\end{equation}
where $f^\ast$ is a non-negative symmetric function which is bounded on $(-\pi,\pi)$ and is bounded away
from zero in a neighborhood of the origin. M($d$) processes encompass both stationary and non-stationary processes, depending
on the value of the memory parameter $d$. When $d < 1/2$, the process $X$ is covariance stationary and its spectral density is given by
\begin{equation}
\label{eq:SpectralDensity:FractionalProcess}
f(\lambda) = |1 - \rme^{-\rmi \lambda}|^{-2d}  f^\ast(\lambda) \eqsp .
\end{equation}
The process $X$ is said to have long-memory if $0 < d < 1/2$, short-memory if $d=0$ and negative memory if $d < 0$; the
process is not invertible if $d < -1/2$. When $d > 1/2$, the process is non stationary. In this case, the $f$
in~(\ref{eq:SpectralDensity:FractionalProcess}) is not integrable on $[-\pi,\pi]$ and is therefore not a spectral density. In
the terminology of~\cite{yaglom:1958}, this $f$ is called a \textit{generalized spectral density}. 
It corresponds to a process $X$ whose increments of sufficiently high order are covariance stationary. 

The memory parameter $d$ plays a central role in the definition of M($d$) processes and is often the focus of empirical interest.
In the parametric case one can use approximate MLE methods (\cite{fox:taqqu:1986}) or MLE (\cite{dahlhaus:1989}). In the
semi-parametric case~(\ref{eq:SpectralDensity:FractionalProcess}) where only a class of functions $f^\ast$ is specified, 
two types of methods have emerged to estimate the memory parameter $d$ : Fourier and Wavelet methods.
Frequency-domain techniques are now well documented and understood  (see for instance
\cite{hurvich:ray:1995}, \cite{velasco:1999}, \cite{velasco:robinson:2000} and \cite{hurvich:moulines:soulier:2002}).

In this paper, we focus on wavelet methods and consider 
the regression estimator introduced in \cite{abry:veitch:1998}, which involves  estimating $d$ using the slope of the regression
of the logarithm of the scale spectrum on the scale index. This estimator is now  widely used in many different fields
(see \eg \cite{veitch:abry:1999} for applications to network traffic; \cite{percival:walden:2000} and \cite{papanicolaou:solna:2003} for applications in physical sciences; see \eg
\cite{gencay:selchuk:whitcher:2002} and \cite{bayraktar:poor:sircar:2004} for applications in finance).
The regression estimator is well-suited to process large data sets, since it has low computational complexity due to the pyramidal algorithm
for computing the details coefficients. Also, it is robust  with respect to additive polynomial trends (see for instance
\cite{veitch:abry:1999} and \cite{craigmile:percival:guttorp:2005}). In Moulines, Roueff and
Taqqu~(\citeyear{moulines:roueff:taqqu:2005b}), we study another estimator of $d$ obtained by adapting the local Whittle
estimator to the wavelet context.

Despite its widespread use, a rigorous semi-parametric asymptotic theory of the regression estimator,
comparable to the one developed for 
corresponding estimators based on the periodogram, is still missing (the concluding remarks 
in \cite{velasco:1999} about ``the lack of rigorous asymptotic theory (...) if the spectral density is not proportional to
$\lambda^{-2d}$ for all frequencies'' for wavelet-based estimates are still valid).   
There are results in a related parametric framework (see \cite{bardet:2002} and \cite{bayraktar:poor:sircar:2004}).
To the best of our knowledge, the only attempt in a semi-parametric setting is due to
\cite{bardet:lang:moulines:soulier:2000}. The process, however, is supposed to be observed in continuous-time
-- discretization issues were not discussed -- and the results do not directly translate to discrete-time observations in a semi-parametric framework. 
The main objective of this paper is to fill this gap. 

The paper is organized as follows. Examples of $M(d)$ processes are given in Section~\ref{sec:examples}. In
Section~\ref{sec:DWT}, we introduce wavelets and wavelet transforms for time-series. 
We do not assume that the wavelets are orthonormal nor that they result from a multiresolution analysis. In
Section~\ref{sec:correlation:wavelet-coefficient}, we derive an explicit expression for the covariance and spectral density
of the wavelet coefficients 
of an M($d$) process   at a given scale.
We extend this result to pairs of scales  by grouping, in an appropriate way, the wavelet coefficients. The results apply to a general class
of wavelets with bounded supports, which include but are not limited to Daubechies wavelets.
If $f^\ast$ belongs to a class of smooth functions $\calH(\beta,L)$ with smoothness exponent $\beta$ defined
in~(\ref{eq:Hbeta}), we show that the spectral density of the wavelet coefficients of  
an M($d$) process can be approximated, at large scales, by the spectral density of the wavelet coefficients of fractional
Brownian motion (FBM)
and derive an explicit bound for the difference between these two quantities. Our result holds not only for $d\in(1/2,3/2)$,
which corresponds to the standard range for the Hurst index, $H= d-1/2\in(0,1)$, 
but for all $d\in\Rset$ by interpreting the corresponding FBM as a generalized process with spectral density
$|\lambda|^{-2d}$, $\lambda\in\Rset$. We show that the relative $L^\infty$ error between the spectral densities of the wavelet
coefficients decreases exponentially fast to zero with a rate given by the smoothness exponent $\beta$ of $f^\ast$. 
In Section \ref{sec:abry-veitch}, we consider (possibly non-stationary) Gaussian processes and obtain
an explicit expression for the limiting variance of the estimator of $d$ based on the regression of the log-scale spectrum.
We show that this estimator is rate optimal in the minimax sense. 
Sections~\ref{sec:proof:prop:VjapproxDensityJbound} and~\ref{sec:proofs-theor-refth} contain
proofs. Appendix~\ref{sec:appr-wavel-filt} involves approximations of wavelet filter transfer functions. We derive 
in Appendix~\ref{sec:useful-inequality} 
an inequality for the mean and the covariance of the logarithm of quadratic
forms of Gaussian variables. 

\section{Examples}\label{sec:examples}

Stationarity of the increments is commonly assumed in time-series analysis. In ARIMA models, for example,~(\ref{eq:spdelta})
holds with $d=K$ integer and with $f^\ast$ equal to the spectral density
of an autoregressive moving average short-memory process.
If $d\in\Rset$ and $f^\ast \equiv \sigma^2$  in~\eqref{eq:SpectralDensity:FractionalProcess},
one gets the so-called fractionally integrated white noise process, ARFIMA(0,d,0). The choice $d\in\Rset$ and
\begin{equation}
\label{eq:smooth:part:ARMA}
f_{\mathrm{ARMA}}^\ast(\lambda)= \sigma^2 \frac{\left| 1 - \sum_{k=1}^q \theta_k \rme^{-\rmi \lambda k}\right|^2}
{\left|1 - \sum_{k=1}^p \phi_k \rme^{- \rmi \lambda k} \right|^2},\quad \lambda\in(-\pi,\pi)\;, 
\end{equation}
with $1 - \sum_{k=1}^p \phi_k z^k \ne 0$ for $|z|=1$ and  $1 - \sum_{k=1}^p \theta_k\ne0$ (so that $f_{\mathrm{ARMA}}^\ast(0)\ne0$) 
leads to the class of ARFIMA($p,d,q$) processes.

Another example is $\{ B_H(k) \}_{k \in \Zset}$,  a discrete-time version of fractional Brownian
motion (FBM) $\{ B_H(t), t \in \Rset \}$ with Hurst index $H \in (0,1)$. The latter is a centered Gaussian process with
covariance 
\[
R_H(t,s) \eqdef  \PE[ B_H(t) B_H(s)] = \frac{1}{2} \left\{ |t|^{2H} + |s|^{2H} - |t-s|^{2H} \right\} \eqsp .
\]
The process $\{ B_H(k) \}_{k \in \Zset}$ is increment stationary ($K=1$) and its generalized spectral density is given up to
a multiplicative constant (see~\cite{samorodnitsky:taqqu:1994}) by
$$
f_{\mathrm{FBM}}(\lambda) \eqdef \sum_{k=-\infty}^\infty\left|\lambda+2k\pi\right|^{-2H-1},\quad \lambda\in(-\pi,\pi)\;.
$$
We can express it in the form~(\ref{eq:SpectralDensity:FractionalProcess}),
\begin{equation}
\label{eq:GenDens:FBM}
f_{\mathrm{FBM}}(\lambda)=|1-\rme^{-\rmi\lambda}|^{-2d}\,f_{\mathrm{FBM}}^\ast(\lambda) \; ,
\end{equation}
by setting $d= H+1/2\in(1/2,3/2)$ and 
\begin{equation}
\label{eq:smooth:part:FBM}
f_{\mathrm{FBM}}^\ast(\lambda)=\left|\frac{2 \sin(\lambda/2)}{\lambda}\right|^{2H+1}+
\left|2\sin(\lambda/2)\right|^{2H+1}\sum_{k\neq0}\left|\lambda+2k\pi\right|^{-2H-1} \; .
\end{equation}
Observe that $f_{\mathrm{FBM}}^\ast(0)=1$ and that it is bounded on $(-\pi,\pi)$.

The process $G_H=\diffop B_H$ is fractional Gaussian noise (FGN). It is a stationary Gaussian process with spectral density
proportional 
to~(\ref{eq:GenDens:FBM}), but with $d=H-1/2\in(-1/2,1/2)$. 

By convention, throughout the paper, while $d$ may take values in $\Rset$, $H$ will be always restricted to take values in
$(0,1)$. 

\section{Discrete Wavelet Transform}
\label{sec:DWT}
In this section, we introduce the main concepts required to define a discrete wavelet transform.
Denote by $\ltwo(\Rset)$ the set of square integrable functions with respect to the Lebesgue measure.
Let $\phi \in \ltwo(\Rset)$ and $\psi \in \ltwo(\Rset)$ be two functions
and define their Fourier transforms as
\[
\hat{\phi}(\xi) \eqdef \int_{-\infty}^\infty \phi(t) \rme^{- \rmi \xi t}\,dt \quad
\text{and}
\quad
\hat{\psi}(\xi) \eqdef \int_{-\infty}^\infty \psi(t) \rme^{- \rmi \xi t}\,dt.
\]
Consider the following assumptions:
\begin{enum_W}
\item\label{item:Wreg} $\phi$ and $\psi$ are compactly-supported, integrable, and $\hat{\phi}(0) = \int_{-\infty}^\infty \phi(t)\,dt = 1$ and $\int_{-\infty}^\infty \psi^2(t)\,dt = 1$.
\item\label{item:psiHat}
There exists $\alpha>1$ such that
$\sup_{\xi\in\Rset}|\hat{\psi}(\xi)|\,(1+|\xi|)^{\alpha} <\infty$.
\item\label{item:MVM} The function $\psi$ has  $M$ vanishing moments, \ie\ $ \int_{-\infty}^\infty t^m \psi(t) \,dt=0$ for all $m=0,\dots,M-1$
\item\label{item:MIM} The function $ \sum_{k\in\Zset} k^m\phi(\cdot-k)$
is a polynomial of degree $m$ for all $m=0,\dots,M-1$.
\end{enum_W}

Assumption \ref{item:Wreg} implies that $\hat{\phi}$ and $\hat{\psi}$ are everywhere infinitely differentiable.
The exponent $\alpha$ in~\ref{item:psiHat} is related to the rate of decrease of the Fourier transform
$\hat{\psi}$ of the wavelet $\psi$ and hence to the regularity of $\psi$.
Under \ref{item:Wreg},~\ref{item:MVM} is equivalent to asserting that the first $M-1$ derivatives of $\hat{\psi}$ vanish at
the origin. This implies, using a Taylor expansion, that
\begin{equation}
\label{eq:MVM}
|\hat{\psi}(\xi)|=O(|\xi|^{M})\quad\text{as}\quad \xi \to 0 \eqsp.
\end{equation}
By \cite[Theorem~2.8.1, Page~90]{cohen:2003}, under~\ref{item:Wreg},~\ref{item:MIM} is equivalent to
\begin{equation}
\label{eq:MIM}
\sup_{k\neq0} |\hat{\phi}(\xi+2k\pi)|=O(|\xi|^{M})\quad\text{as}\quad\xi\to0 \eqsp.
\end{equation}

Adopting the engineering convention that large values of the scale index $j$ correspond to coarse scales (low frequencies),
we define the family $\{\psi_{j,k}, j \in \Zset, k \in \Zset\}$ of translated and dilated functions 
\begin{equation}\label{eq:psiJK}
\psi_{j,k}(t)=2^{-j/2}\,\psi(2^{-j}t-k) \eqsp .
\end{equation}
Many authors suppose that the $\psi_{j,k}$ are orthogonal and even that they are generated by a multiresolution analysis (MRA).
Assumptions \allWA\ are indeed quite standard in the context of a multiresolution analysis (in which case, $\phi$ is the 
scaling function and $\psi$ is the associated wavelet), see for instance~\cite{cohen:2003}. 
In this paper, we do not assume that wavelets are orthonormal nor
that they are associated to a multiresolution analysis. We may therefore work with other convenient choices for $\phi$ and $\psi$
as long as~\allWA\ are satisfied.  A simple example is to set, for some positive integer $N$,
$$
\phi(x) \eqdef \1_{[0,1]}^{\otimes N}(x)\quad\text{and}\quad\psi(x) \eqdef \frac{d^N}{d x^N} \1_{[0,1]}^{\otimes 2N}(x),
$$
where $\1_A$ is the indicator function of the set $A$ and for a non negative function $f$, $f^{\otimes N}$ denotes the $N$-th  self-convolution of
$f$. It follows that
$$
|\hat\phi(\xi)|=|2\sin(\xi/2)/\xi|^{N}\quad\text{and}\quad\hat\psi(\xi)=|\xi|^{N}|2\sin(\xi/2)/\xi|^{2N}.
$$
Using~(\ref{eq:MVM}) and~(\ref{eq:MIM}), one easily checks that~\allWA\ are satisfied with $M$ and $\alpha$ equal to $N$.
Of course the  family of functions $\{\psi_{j,k}\}$ are not orthonormal for this choice of the wavelet function $\psi$ (and the function
$\phi$ is not associated to a MRA). Nevertheless, to ease references to previously reported works,  with a slight abuse in the terminology,
we still call  $\phi$ and $\psi$ the scaling and the wavelet functions.

Having defined the functions $\phi$ and $\psi$, we now define what we call the Discrete Wavelet Transform in discrete
time. Start with a real-valued sequence $\{x_k,\,k\in\Zset\}$. Using the scaling function $\phi$, we first associate to the
sequence $\{x_k, k \in \Zset\}$ the continuous-time functions
\begin{equation}\label{eq:bX}
\bx_n(t) \eqdef \sum_{k=1}^n x_k \,\phi(t-k) \quad\text{and}\quad \bx(t) \eqdef \sum_{k\in\Zset} x_k\, \phi(t-k), \quad
t\in\Rset \; .
\end{equation}
The wavelet coefficients involve $\{\bx(t),\,t\in\Rset\}$ and are defined as
\begin{equation}\label{eq:coeff}
\dwt^\bx_{j,k} \eqdef \int_{-\infty}^\infty \bx(t) \psi_{j,k}(t)\,dt\quad  j \geq 0, k \in \Zset.
\end{equation}

Without loss of generality we may suppose that the support of the scaling function $\phi$ is included in $(-\L,0)$ for some
integer $\L\geq1$. Then $\bx_n(t)=\bx(t)$ for all $t\in(0, n-\L+1)$. 
We may also suppose that the support of the wavelet function $\psi$ is included in $(0,\L)$. Then, the support of $\psi_{j,k}$ is included in the
interval $(2^j k, 2^j(k+\L) )$. Hence
\begin{equation}\label{eq:coeffN}
\dwt^\bx_{j,k}=\dwt^{\bx_n}_{j,k}=\int_{-\infty}^\infty \bx_n(t) \psi_{j,k}(t)\,dt,
\end{equation}
when $(2^j k, 2^j(k+\L) )\subseteq (0, n-\L+1)$, that is, for all $(j,k)\in\indexset_n$, where
\begin{equation}
\label{eq:deltan}
\indexset_n \eqdef \{(j,k):\,j\geq0, 0\leq k  \leq 2^{-j}(n-\L+1)-\L \} \eqsp.
\end{equation}
For any $j$, the wavelet coefficients $\{\dwt^\bx_{j,k} \}_{k \in \Zset}$ are obtained by discrete convolution and downsampling.
More precisely, under \ref{item:Wreg}, for all $j \geq 0$, $k \in \Zset$,~(\ref{eq:bX}) and~(\ref{eq:coeff}) imply
\begin{equation}\label{eq:down}
\dwt^{\bx}_{j,k}=\sum_{l\in\Zset} x_l\,h_{j,2^jk-l}=(h_{j,\cdot}\star \bx)_{2^jk}= (\downarrow^j [h_{j,\cdot}\star \bx])_{k},
\end{equation}
where $\star$ denotes the convolution  of discrete sequences,
$$
 h_{j,l} \eqdef 2^{-j/2} \int_{-\infty}^\infty \phi(t+l)\psi(2^{-j}t)\,dt \; ,
$$  
and, for any sequence $\{ c_k \}_{k \in \Zset}$ and any integer $l$,
$\downarrow^j$ is the downsampling operator defined as  
$(\downarrow^j c)_k = c_{2^j k}$.
Define, for all $j\geq0$, the discrete Fourier transform of $\{ h_{j,l} \}_{l \in \Zset}$ as
\begin{equation}\label{eq:fjdef}
H_j(\lambda) \eqdef \sum_{l\in\Zset}h_{j,l}\rme^{-\rmi\lambda l} =2^{-j/2} \int_{-\infty}^\infty
\sum_{l\in\Zset}\phi(t+l)\rme^{-\rmi\lambda l} \psi(2^{-j}t)\,dt. 
\end{equation}
Since $\phi$ and $\psi$ have compact support, the sum in~(\ref{eq:fjdef}) has a finite number of non-vanishing terms and
$H_j$ is a trigonometric polynomial. 

\begin{remark}\label{rem:j0}
By Corollary~\ref{cor:InfHj}, there exists an integer $j_{0}$ only depending on $\phi$
and $\psi$ such that, for all $j\geq j_{0}$, the trigonometric polynomial $H_j(\lambda)$ is not identically zero. 
In the case of a multiresolution analysis, the   
father and mother wavelets are defined in such a way that $j_0=0$. 
In the general case, by dilating $\psi$ appropriately, or, in other words, by changing the reference scale,
one can impose $j_0=0$, which is assumed in the sequel.    
\end{remark} 

Under assumption \ref{item:MIM}, $t \mapsto \sum_{l \in \Zset} \phi(t+l) l^m$ is a polynomial of degree $m$ and \ref{item:MVM} therefore
implies that, for all $j\geq0$ and all $m=0,\dots,M-1$, 
\begin{equation}
\label{eq:polyHj}
\sum_{l\in\Zset} h_{j,l}\,l^m= 2^{-j/2} \int_{-\infty}^\infty \psi(2^{-j}t) \sum_{l\in \Zset} \phi(t+l) l^m dt =0 \eqsp.
\end{equation}
Now consider $P_j(x)=\sum_{l \in \Zset} h_{j,l} \, x^l$ and observe that~(\ref{eq:polyHj}) implies $P_j(1)=0$, $P_j'(1)=0$,
..., $P_j^{(M-1)}(1)=0$, and hence $H_j(\lambda)=P_j(\rme^{-\rmi\lambda})$ factors as
\begin{equation}
\label{eq:Htilde}
 H_{j}(\lambda)= (1-\rme^{-\rmi \lambda})^{M} \, \tilde{H}_j(\lambda) \; ,
\end{equation}
where $\tilde{H}_j(\lambda)$ is also a trigonometric polynomial. The wavelet coefficient \eqref{eq:down} may therefore be computed as
\begin{equation}
\label{eq:op}
\dwt^{\bx}_{j,k}= (\downarrow^j [\tilde{h}_{j,\cdot}\star \diffop^M \bx])_{k}
\end{equation}
where $ \{ \tilde{h}_{j,l} \}_{l \in \Zset}$ are the coefficients of the trigonometric polynomial $\tilde{H}_j$ and
$\diffop^M \bx$ is the $M$-th order difference of the sequence $\bx$.
In other words, the use of a wavelet and a scaling function satisfying \ref{item:MIM} and \ref{item:MVM} implicitly perform a $M$-th order differentiation
of the time-series. Therefore, we may work with a $K$-th order integrated processes $X$ without specific preprocessing, provided that $M \geq K$. This
is in sharp contrast with Fourier methods, where the time series must be explicitly differentiated at least $K$ times and a
data taper must be applied on the differenced series to avoid frequency-domain leakage
(see, for instance, \cite{hurvich:moulines:soulier:2002}).


\section{Spectral Density of the Wavelet Coefficients}
\label{sec:correlation:wavelet-coefficient}
Because the wavelet coefficients at a given scale are obtained by applying time-invariant linear filters, computing the covariance of
the wavelet coefficients of $K$-th order stationary processes is an easy exercise.
The following proposition provides an integral expression for the covariance between two
wavelet coefficients on possibly different scales, expressed in terms of the transfer function $H_j$ of the
wavelet filters and the generalized spectral density of the process $X$. This proposition extends
Theorem~2 in \cite{masry:1993} on the spectral measure of the DWT coefficients of increment stationary continuous-time
processes  to the discrete-time setting and Lemma~1 in \cite{craigmile:percival:2005} to functions $\psi$ and $\phi$ that do
not necessarily define a multiresolution analysis. 
\begin{proposition}\label{prop:covdjk}
Let $X$ be a $K$-th order integrated process with generalized
spectral density $f$. Assume \allWA\ with $M\geq K$. Then, for all
$j,j' \geq 0$ and $k,k' \in \Zset$,
\begin{equation}\label{eq:covWC}
\PCov(\dwt^X_{j,k},\dwt^X_{j',k'})= \int_{-\pi}^\pi \rme^{\rmi \lambda (k2^j-k'2^{j'})}\,f(\lambda)\,
H_j(\lambda)\,\overline{H_{j'}(\lambda)}\,d\lambda,
\end{equation}
where the wavelet coefficient $\dwt_{j,k}^X$ is defined in \eqref{eq:coeff}.
\end{proposition}
The proof follows from elementary results on time-invariant linear filtering of covariance stationary processes,
using~(\ref{eq:Htilde}),~(\ref{eq:op}), applied to $\diffop^MX$, which is covariance stationary with spectral density
$|1-\rme^{\rmi \lambda}|^{2M}f(\lambda)$.

By \eqref{eq:op}, for a given scale $j$, the process $\{ \dwt^X_{j,k} \}_{k \in \Zset}$ is covariance stationary.
The situation is more complicated when considering two different scales $j \ne j'$, because
the two-dimensional sequence $\{[\dwt^X_{j,k},\,\dwt^X_{j',k}]^T \}_{k\in\Zset}$, with $^T$ denoting the transpose, is not
stationary for $j\neq j'$. 
This is a consequence of the pyramidal wavelet scheme, where at scale $j$,
the wavelet coefficients are downsampled by a factor $2^j$ which depends on $j$ since $\psi_{j,k}$ in~(\ref{eq:psiJK}) can be
expressed as 
$$
\psi_{j,k}(t)=2^{-j/2}\psi(2^{-j}(t-2^jk))\;.
$$
Thus, to obtain a stationary sequence, one should
consider the process
$\{[\dwt^X_{j,k},\,\dwt^X_{j',2^{j-j'}k}]^T \}_{k\in\Zset}$ for $j> j'$, which involves a downsampled subsequence of the
coefficients at the finer scale $j'$. One can also consider the process
$\{[\dwt^X_{j,k},\,\dwt^X_{j',2^{j-j'}k+\dk}]^T \}_{k\in\Zset}$ for $j> j'$, which includes a
translation of the location index of the second component by $\dk$. It turns out that the most convenient is to
merge the processes corresponding to $\dk=0,\dots,2^{j-j'}-1$ and hence to consider the \emph{between-scale} process
\begin{equation}\label{eq:betweenscaleProc}
\{[\dwt^X_{j,k},\,\bdwt^X_{j,k}(j-j')^T]^T \}_{k\in \Zset}\;, 
\end{equation}
where for any $\dj=0,1,\dots,j$,
\begin{equation}\label{eq:Defbd}
\bdwt^X_{j,k}(\dj) \eqdef \left[\dwt^X_{j-\dj,2^{\dj}k},\,\dwt^X_{j-\dj,2^{\dj}k+1},\,\dots, \dwt^X_{j-\dj,2^{\dj}k+2^{\dj}-1}\right]^T 
\end{equation}
is a $2^\dj$-dimensional vector of wavelet coefficients at scale $j'=j-\dj$. The vector $\bdwt^X_{j,k}(\dj)$ involves all
possible translations of the position index $2^{\dj}k$ by $\dk=0,1,\dots,2^{\dj}-1$.
The index $\dj$ in~(\ref{eq:Defbd}) denotes the scale difference $j-j'\geq0$ between the finest scale $j'$ and the coarsest
scale $j$. Observe that $\bdwt^X_{j,k}(0)$ ($\dj=0$) is the scalar $\dwt^X_{j,k}$.

One should view the between-scale process~(\ref{eq:betweenscaleProc}) as a pair made up of the scalar process
$\{\dwt^X_{j,k}\}_{k\in\Zset}$ and the vector process $\{\bdwt^X_{j,k}(j-j')\}_{k \in\Zset}$. We shall now express their
cross spectral density in 
terms of the generalized spectral density of $X$ and the transfer function of the wavelet filters folded on the interval
$[-\pi,\pi]$. By setting $j'=j$ or equivalently $\dj=0$ we obtain the spectral density of the ``within scale'' process
$\{\dwt^X_{j,k}\}_{k\in\Zset}$. 
\begin{corollary} \label{cor:bD}
Define for all $0 \leq \dj \leq j$ and $\lambda \in (-\pi,\pi)$,
\begin{multline}
\label{eq:definitionfj}
\bdens[\phi,\psi]{j,\dj}{\lambda}{f} \eqdef  \\ 
\sum_{l=0}^{2^{j}-1}  \be_{\dj}(\lambda+2l\pi)\, f(2^{-j}(\lambda+2l\pi))\,
2^{-j/2} H_{j}(2^{-j}(\lambda+2l\pi))\, 2^{-(j-\dj)/2}\overline{H_{j-\dj}(2^{-j}(\lambda+2l\pi))} \eqsp,
\end{multline}
where for all $\xi\in\Rset$, 
$$ 
\be_\dj(\xi) \eqdef 2^{-\dj/2}\, [1, \rme^{-\rmi2^{-\dj}\xi}, \dots, \rme^{-\rmi(2^{\dj}-1)2^{-\dj}\xi}]^T\;.
$$
Then
$$
\PCov_f(\dwt^X_{j,k},\bdwt^X_{j,k'}(\dj)) = \int_{-\pi}^\pi \rme^{\rmi\lambda(k-k')} \,
\bdens[\phi,\psi]{j,\dj}{\lambda}{f} \, d\lambda \; .
$$
In other words, 
\begin{itemize}
\item for all $j \geq 0$, the within-scale process $ \{ \dwt^X_{j,k} \}_{k \in \Zset}$ is covariance stationary with spectral density
$\bdens[\phi,\psi]{j,0}{\cdot}{f}$,
\item for all $j \geq \dj > 0$, the between-scale process $\{[\dwt^X_{j,k},\,\bdwt^X_{j,k}(\dj)^T]^T\}_{k\in \Zset}$ is covariance stationary with 
cross spectral density $\bdens[\phi,\psi]{j,\dj}{\cdot}{f}$.
\end{itemize}
\end{corollary}
\noindent Note that $\densletter_{j,\dj}$ is a $2^{\dj}$-dimensional vector and, in particular, $\densletter_{j,0}$ is a
scalar. The $2^\dj$-dimensional vector $\be_\dj(\xi)$ has Euclidean norm $|\be^\dj(\xi)|=1$. 
\begin{proof}
Let  $j\geq \dj\geq0$.
By~(\ref{eq:covWC}), we have, for all $k,k'\in\Zset$ and $\dk=0,\dots,2^{\dj}-1$,
$$
\PCov_f(\dwt^X_{j,k},\dwt^X_{j-\dj,k'2^{\dj}+\dk})=
\int_{0}^{2\pi} \rme^{\rmi\lambda(k2^j- (k'2^{\dj}+\dk)2^{j-\dj})}\,
f(\lambda)\,H_j(\lambda)\,\overline{H_{j-\dj}(\lambda)}\,d\lambda.
$$
The exponential can be factorized as $\rme^{\rmi2^j\lambda(k-k')}\,\rme^{-\rmi\dk2^{-\dj}2^{j}\lambda}$. Hence
\begin{multline*}
\PCov_f(\dwt^X_{j,k},\bdwt^X_{j,k'}(\dj))  =
\int_{0}^{2\pi} \rme^{\rmi\lambda2^{j}(k-k')}\,\be_{\dj}(2^{j}\lambda) \,f(\lambda)\,
H_{j}(\lambda)2^{\dj/2}\overline{H_{j-\dj}(\lambda)}\,d\lambda \\
=
\int_{0}^{2^{j+1}\pi} \rme^{\rmi\lambda(k-k')}\,\be_{\dj}(\lambda) \,f(2^{-j}\lambda)\,
2^{-j/2}H_{j}(2^{-j}\lambda)\,2^{-(j-\dj)/2}\overline{H_{j-\dj}(2^{-j}\lambda)}\,d\lambda.
\end{multline*}
The result is obtained by \textit{folding} and shifting the previous integral as follows, 
setting $g(\lambda)= \be_{\dj}(\lambda) 
\,f(2^{-j}\lambda)\,2^{-j/2}H_{j}(2^{-j}\lambda)\,2^{-(j-\dj)/2}\overline{H_{j-\dj}(2^{-j}\lambda)}$,  
\begin{multline*}
\int_{0}^{2^{j+1}\pi} \rme^{\rmi\lambda(k-k')}\, g(\lambda)\,d\lambda
=\sum_{l=0}^{2^{j}-1} \int_{2l\pi}^{2(l+1)\pi} \rme^{\rmi\lambda(k-k')}\, g(\lambda)\,d\lambda\\
=\sum_{l=0}^{2^{j}-1} \int_{0}^{2\pi} \rme^{\rmi\lambda(k-k')}\,g(\lambda+2l\pi)\,d\lambda
= \int_{0}^{2\pi} \rme^{\rmi\lambda(k-k')}\, \left(\sum_{l=0}^{2^{j}-1}  g(\lambda+2l\pi)\right)\,d\lambda\; .
\end{multline*}
The function in parentheses is $(2\pi)$-periodic because 
$\sum_{l=0}^{2^{j}-1}  g(\lambda+2(l+1)\pi)=\sum_{l=1}^{2^{j}-1}  g(\lambda+2l\pi)+ g(\lambda+2^j(2\pi))=
 \sum_{l=0}^{2^{j}-1}  g(\lambda+2l\pi)$ since $g$ is $2^{j}(2\pi)$-periodic. Hence $\int_{0}^{2\pi}$ can be replaced by
 $\int_{-\pi}^{\pi}$, which gives the result. 
\end{proof}

We now apply the preceding results to the class of processes with memory parameter $d \in \Rset$ (see \eqref{eq:spdelta}).
To obtain error bounds on the variance and the spectral density of the wavelet coefficients, some additional
assumptions are required on the smoothness of $f^\ast$ at zero frequency. 
\begin{definition}
For $0 < \beta \leq 2$ and $L>0$, define the function class
$\calH(\beta,L)$ as the set of even non-negative functions $g$ on  $[-\pi,\pi]$  such that, for all 
$\lambda\in[-\pi,\pi]$,
\begin{equation}\label{eq:Hbeta}
|g(\lambda)- g(0) | \leq L \, g(0) \, |\lambda|^\beta \eqsp .
\end{equation}
\end{definition}
%
This type of assumption is  typical in the \textit{semi-parametric} estimation setting (see for instance
\cite{robinson:1995g} and \cite{moulines:soulier:2001}). The larger the value of $\beta$, the smoother the function at the
origin. For $g$ even -- as assumed -- and infinitely differentiable, $g'(0)=0$ and hence, by a Taylor
expansion,~(\ref{eq:Hbeta}) holds with $\beta=2$.

Since, for instance, $f^\ast_{\mathrm{ARMA}}$ in~\eqref{eq:smooth:part:ARMA} is infinitely
differentiable, it belongs to $\calH(2,L)$ for some $L$. 
As for FBM and FGN, consider $f_{\mathrm{FBM}}^\ast$ in~\eqref{eq:smooth:part:FBM}. The first term in the RHS is
$1+O(\lambda^2)$ and the second is $O(|\lambda|^{2H+1})$. Hence, for some positive constant $L$,
$|f_{\mathrm{FBM}}^\ast(\lambda)-f_{\mathrm{FBM}}^\ast(0)|\leq L\,f_{\mathrm{FBM}}^\ast(0)\,|\lambda|^{(2H+1)\wedge2}$, where
$a\wedge b=\min(a,b)$; hence 
\begin{equation}
\label{eq:fstarFBM:calh}
f^\ast_{\mathrm{FBM}}\in\calH((2H+1) \wedge 2,L)\;.
\end{equation}

The expressions of the within- and between-scale wavelet coefficient  spectral densities
$\bdens[\phi,\psi]{j,\dj}{\cdot}{f}$  given in Corollary~\ref{cor:bD}
depend both on $d$ and on the function $f^\ast$ and will therefore be denoted by $\bdens[\phi,\psi]{j,\dj}{\cdot}{d,f^\ast}$ in
the sequel. We are going
to show, however, that these quantities may be approximated by quantities which depend only on the memory parameter $d$ and $f^\ast(0)$.
Let $X$ have a generalized spectral density $f(\lambda) = |1- \rme^{\rmi \lambda}|^{-2d} f^\ast(\lambda)$ and define
\begin{equation}
\label{eq:definitionvj}
\vj{j}{d,f^{\ast}} \eqdef \PVar [ \dwt^X_{j,0}] = \int_{-\pi}^{\pi} |1 - \rme^{-\rmi \lambda}|^{-2 d} f^\ast(\lambda) |H_j(\lambda)|^2  d \lambda \eqsp,
\end{equation}
the variance of the wavelet coefficient of the process $X$ at scale $j$.
\begin{theorem}\label{theo:VjapproxDensityJbound}
Let $M \geq 1$ be an integer and  $\alpha$, $L$, $\beta$ be constants such that $\alpha > 1$, $ 0 < L < \infty$ and $\beta \in (0,2]$.
Assume that \allWA\ hold with $M$ and $\alpha$.
\begin{enum_a}
\item \label{item:VjApprox} Let $\dmin$ and $\dmax$ be two constants such that
\begin{equation}\label{eq:CondDminDmax-a}
[\dmin, \dmax] \subset \left( (1+\beta)/2-\alpha , M+1/2 \right) \eqsp.
\end{equation}
Then, there exists a constant $C>0$ (only depending on the constants $\beta$, $\dmin$, $\dmax$ and the functions $\phi$ and $\psi$) such that,
for all $j\geq 0$, $d \in [\dmin,\dmax]$ and $f^\ast \in \calH(\beta,L)$,
\begin{equation}\label{eq:VjApprox}
\left| \vj{j}{d,f^\ast} - f^\ast(0) \, \Kvar[\psi]{d} \, 2^{2jd} \right| \leq C\, f^\ast(0) \, L \, 2^{(2d-\beta)j}
\end{equation}
where $\Kvar[\psi]{d}$ is given by
\begin{equation}\label{eq:Kpsi}
\Kvar[\psi]{d}\eqdef \int_{-\infty}^{\infty} |\xi|^{-2d}\,|\hat\psi(\xi)|^2\,d\xi \eqsp.
\end{equation}
\item \label{item:DjApprox} Let $\dmin$ and $\dmax$ be two constants such that
\begin{equation}\label{eq:CondDminDmax-b}
[\dmin, \dmax] \subset \left( (1+\beta)/2-\alpha , M \right] \eqsp.
\end{equation}
Then, for all $\dj\geq0$, there exists $C>0$ (only depending on $\dj$ and on the constants $\beta$, $\dmin$, $\dmax$ and the
functions $\phi$ and $\psi$) such that, for all $\lambda\in(-\pi,\pi)$, $j\geq \dj$, $d \in
[\dmin, \dmax]$ and  $f^\ast \in \calH(\beta,L)$, 
\begin{equation}\label{eq:DjApprox}
\left|\bdens[\phi,\psi]{j,\dj}{\lambda}{d,f^\ast} - f^\ast(0)\, \bdensasymp[\psi]{\dj}{\lambda}{d} \, 2^{2j d}\right|
\leq C\, f^\ast(0) \, L \, 2^{(2d-\beta)j}
\end{equation}
where $|\cdot|$ denotes the Euclidean norm in any dimension and, for all $\dj\geq0$,
\begin{equation}\label{eq:bDpsi}
\bdensasymp[\psi]{\dj}{\lambda}{d} \eqdef  
\sum_{l\in\Zset} |\lambda+2l\pi|^{-2d}\,\be_{\dj}(\lambda+2l\pi) \,
\overline{\hat{\psi}(\lambda+2l\pi)}\hat{\psi}(2^{-\dj}(\lambda+2l\pi)).
\end{equation}
The function $(\lambda,d)\mapsto\bdensasymp[\psi]{\dj}{\lambda}{d}$ is $(2\pi)$-periodic in $\lambda$ and jointly
continuous on $\Rset\times [\dmin,\dmax]$. When $\dj=0$, $\densletter_{\infty,0}$ is a scalar and 
\begin{equation}\label{eq:bDinftyDoNotVanish}
\int_{-\pi}^\pi
\bdensasymp[\psi]{0}{\lambda}{d}\,d\lambda =\Kvar[\psi]{d}\neq0 \; .
\end{equation}
\end{enum_a}
\end{theorem}
The proof, based on approximating the wavelet filter transfer function, can be found in
Section~\ref{sec:proof:prop:VjapproxDensityJbound}.
In order to shed light on Theorem~\ref{theo:VjapproxDensityJbound}, we conclude this section with a number of remarks.
\begin{remark}
Theorem~\ref{theo:VjapproxDensityJbound} states that $f^\ast(0) \, \Kvar[\psi]{d} \, 2^{2jd}$ is a good approximation for 
$\PVar [ \dwt^X_{j,0}]$ and that for any $\dj \geq 0$,
$f^\ast(0) \bdensasymp[\psi]{\dj}{\lambda}{d} 2^{2jd}$ is a good $L^\infty(-\pi,\pi)$ approximation to the spectral density
$\bdens[\phi,\psi]{j,\dj}{\lambda}{d,f^\ast}$.
\end{remark}
\begin{remark}
Relation~(\ref{eq:DjApprox}) with $\dj=0$ implies~(\ref{eq:VjApprox}) since
$\left|\int_{-\pi}^\pi g_1(\lambda)d\lambda-\int_{-\pi}^\pi g_2(\lambda)d\lambda\right|
\leq 2\pi\|g_1-g_2\|_\infty$. Observe, however, that~(\ref{eq:VjApprox}) is valid under Condition~(\ref{eq:CondDminDmax-a}),
which is weaker than~(\ref{eq:CondDminDmax-b}).
\end{remark}
\begin{remark}\label{rem:BdAsympLpNorm}
Under Condition~(\ref{eq:CondDminDmax-b}), for all $p>0$,
\begin{equation}\label{eq:BdAsympLpNorm}
\|\bdensasymp[\psi]{d}{\cdot}{d}\|_p
\eqdef
\left(\int_{-\pi}^\pi \left|\bdensasymp[\psi]{d}{\lambda}{d}\right|^p\,d\lambda \right)^{1/p}
\end{equation}
is positive, finite and continuous in $d\in[\dmin,\dmax]$. This follows from joint continuity and~(\ref{eq:bDinftyDoNotVanish}). 
\end{remark}
\begin{remark}
The spectral density $f^\ast(0)\, \bdensasymp[\psi]{\dj}{\lambda}{d} \, 2^{2j d}$, $\lambda\in(-\pi,\pi)$, is in fact the
spectral density of the wavelet coefficient of the generalized fractional Brownian motion $B_{(d)}$, where
$d\in\Rset$. The process $B_{(d)}$ is parameterized by a family $\Theta_{(d)}$ of smooth ``test'' functions $\theta(t)$,
$t\in\Rset$ and is defined as follows: $\{B_{(d)}(\theta),\,\theta\in\Theta_{(d)}\}$ is a mean zero Gaussian process
with covariance
\begin{equation}\label{eq:CovFBMGen}
\PCov\left(B_{(d)}(\theta_1),\,B_{(d)}(\theta_2)\right) = \int_{\Rset}
|\lambda|^{-2d}\,\hat{\theta_1}(\lambda)\,\overline{\hat{\theta_2}(\lambda)}\,d\lambda\;. 
\end{equation}
The finiteness of the integral $\int_\Rset|\lambda|^{-2d}\,|\hat{\theta}(\lambda)|^2\,d\lambda$ provides a constraint on the
family $\Theta_{(d)}$. For instance, when $d>1/2$, this condition requires that $\hat{\theta}(\lambda)$ decays sufficiently
quickly at the origin and, when $d<0$, it requires that $\hat{\theta}(\lambda)$ decreases sufficiently rapidly at
infinity. Hence $\theta$ can be a wavelet $\psi$ if $d\in(1/2-\alpha,M+1/2)$, which corresponds
to~(\ref{eq:CondDminDmax-a}) with $\beta=0$. The discrete wavelet transform of $B_{(d)}$ is defined as
\begin{equation}\label{eq:DWTFBMGen}
\dwt_{j,k}^{(d)} \eqdef B_{(d)}(\psi_{j,k}),\quad j\in\Zset,\,k\in\Zset\;. 
\end{equation}
The spectral density $f^\ast(0) \bdensasymp[\psi]{\dj}{\lambda}{d} 2^{2jd}$ in~(\ref{eq:DjApprox}) which serves as an
approximation to $\bdens[\phi,\psi]{j,\dj}{\lambda}{d,f^\ast}$ is in fact, up to the multiplicative constant $f^\ast(0)$, the cross spectral
density between the wavelet coefficients $\dwt_{j,k}^{(d)}$ and the vector of wavelet coefficients
$\bdwt^{(d)}_{j,k}(\dj)\eqdef\left[\dwt^{(d)}_{j-\dj,2^{\dj}k},\dots, \dwt^{(d)}_{j-\dj,2^{\dj}k+2^{\dj}-1}\right]$. Indeed,
using~(\ref{eq:CovFBMGen}),~(\ref{eq:DWTFBMGen}) and $\hat{\psi}_{j,k}(\lambda)=2^{j/2}\hat{\psi}(2^j\lambda)\rme^{-\rmi
  k\lambda2^j}$, one has 
\begin{align}
\nonumber
\PCov\left(\dwt_{j,k}^{(d)},\,\bdwt^{(d)}_{j,k'}(\dj)\right)& =2^{j} \, \int_\Rset
|\lambda|^{-2d}\,\be_\dj(-2^j\lambda)\,\hat{\psi}(2^{j}\lambda)\overline{\hat{\psi}(2^{j-\dj}\lambda)}\,
\rme^{\rmi \lambda 2^j(k'-k)}\,d\lambda \\
\label{eq:DasympAndCov}
& = 2^{2dj}\,\int_{-\pi}^\pi \bdensasymp[\psi]{\dj}{\lambda}{d} \, \rme^{\rmi \lambda (k-k')}\,d\lambda \; ,
\end{align}
where the last equality is obtained by the change of variable $-2^j\lambda\to\lambda$ and by folding the integral on $(-\pi,\pi)$.
 
The within- and between-scale spectral densities $\bdens[\phi,\psi]{j,\dj}{\lambda}{d,f^\ast}$  
of the process $X$ with memory parameter $d$ may thus be approximated by the DWT of the generalized FBM $B_{(d)}$ 
(viewed as a generalized process), with an $L^\infty$ error bounded by the RHS in~(\ref{eq:DjApprox}).   
\end{remark}

\begin{remark}
When $d$ belongs to $(1/2,3/2)$, $B_{(d)}$ is related to $B_H(t)$, $t\in\Rset$, by setting
$H=d-1/2\in(0,1)$ and, up to a multiplicative constant, 
$$
B_{(d)}(\theta) = \int_{\Rset} B_H(t) \, \theta(t)\, dt \;,
$$
where the equality holds in the sense of finite-dimensional distributions and hence 
$\{\dwt^{(d)}_{j,k},\, j,k\in\Zset\}$ has same distribution as $\{\int_{-\infty}^\infty B_H(s) \psi_{j,k}(s)\,ds,\, j,k\in\Zset\}$.
It follows from the previous remark that, for such $d$ and $H$, the spectral density of the wavelet coefficients of $X$ can be approximated by that
of continuous-time fractional Brownian motion.
\end{remark}
\begin{remark}
Once normalized by $2^{2jd}$, which is the order of the variance of the wavelet coefficients at scale $j$
(see~(\ref{eq:VjApprox})), the difference of the spectral densities in~(\ref{eq:DjApprox}) is bounded by a constant times
$2^{-\beta j}$, a factor which tends to zero exponentially fast as $j\to\infty$. The rate of the decrease is determined by
the smoothness exponent $\beta$ of $f^\ast$. 
\end{remark} 
\begin{remark}
If $d=0$ and $\{\psi_{j,k},\,k\in\Zset, j \in \Zset\}$ is an
orthonormal system, then~(\ref{eq:CovFBMGen}) and~(\ref{eq:DWTFBMGen}) imply
$$
\PCov(\dwt_{j,k}^{(0)},\dwt_{j',k'}^{(0)})=
\int_{\Rset}\hat{\psi}_{j,k}(\lambda)\,\overline{\hat{\psi}_{j',k'}(\lambda)}\,d\lambda
=2\pi\,\int_{\Rset}\psi_{j,k}(t)\psi_{j',k'}(t)\,dt 
$$ 
which vanishes if $j\ne j'$ or $k\ne k'$. 
Hence, when the memory parameter $d=0$ and the wavelets are orthonormal,
the wavelet coefficients $\{\dwt^X_{j,k},\,k\in\Zset\}$ are then asymptotically uncorrelated as $j\to\infty$ and their asymptotic
variance is $2\pi$. Using~(\ref{eq:DasympAndCov}), the corresponding cross spectral density is given by 
\begin{equation}\label{eq:bdensasympNullD}
\bdensasymp[\psi]{\dj}{\lambda}{0}=0\quad\text{if}\quad\dj>0\quad\text{and}\quad
\bdensasymp[\psi]{0}{\lambda}{0}=1,\quad \lambda\in(-\pi,\pi) \; .
\end{equation}
\end{remark}

\begin{remark}
To understand the presence of the asymptotic form $f^\ast(0)2^{2jd}\bdensasymp[\psi]{\dj}{\lambda}{d}$ in~(\ref{eq:DjApprox}), start with
$\bdens[\phi,\psi]{j,\dj}{\lambda}{f}$ in~(\ref{eq:definitionfj}), use the $2^{j}(2\pi)$-periodicity to replace 
$\sum_{l=0}^{2^j-1}$ by $\sum_{l=-2^{j-1}+1}^{2^{j-1}}$, and as $j\to\infty$,
approximate $f$ by the spectral density 
$f^\ast(0)|\lambda|^{-2d}$ of $\sqrt{f^\ast(0)}B_{(d)}$, $H_j(\lambda)$ by its asymptotic approximation 
$2^{j/2}\hat{\phi}(\lambda)\hat{\psi}(2^j\lambda)$ in~(\ref{eq:fjApprox}) and  $H_{j-\dj}(\lambda)$ by
$2^{(j-\dj)/2}\hat{\phi}(\lambda)\hat{\psi}(2^{j-\dj}\lambda)$, approximate $\hat{\phi}(2^{-j}(\lambda+2l\pi))$ by
$\hat{\phi}(0)=1$ and approximate $\sum_{l=-2^{j-1}+1}^{2^{j-1}}$ by $\sum_{l=-\infty}^\infty$.  
\end{remark}

\begin{remark}
Let us examine how Theorem~\ref{theo:VjapproxDensityJbound} applies when $X(k)=B_H(k)$, $k\in\Zset$, that is, $X$ is a
discrete-time version of FBM with Hurst index $H\in(0,1)$. 
From~(\ref{eq:GenDens:FBM}),~(\ref{eq:smooth:part:FBM}) and~(\ref{eq:fstarFBM:calh}), we have $d=H+1/2\in(1/2,3/2)$ and
$f^\ast\in\calH((2H+1)\wedge 2,L)$ for some constant $L$. The condition on $M$ is then $M>H$ in case~\ref{item:VjApprox} and
$M\geq H+1/2$ in case~\ref{item:DjApprox}. The condition on $\alpha$ is $\alpha>\beta/2-H=(1-H)\wedge1/2$
in both cases, which is satisfied because $\alpha>1$ and
$H\in(0,1)$. Theorem~\ref{theo:VjapproxDensityJbound} 
can therefore be applied irrespectively of the value of $H$ when $\psi$  is a Daubechies wavelet with at least $M=2$
vanishing moments. 
\end{remark}   

\begin{remark}
In case of FGN, $d=H-1/2\in(-1/2,1/2)$ and hence, compared to the previous case, $M$ decreases by 1 and $\alpha$
increases by 1. Thus the conditions are $M>H-1$ in case~\ref{item:VjApprox} and
$M\geq H-1/2$ in case~\ref{item:DjApprox}, that is $M=1$ works in either case. The condition on $\alpha$ becomes
$\alpha>(2-H)\wedge3/2$, that is $\alpha>3/2$ will work for all $H\in(0,1)$. Since the Daubechies wavelet
with $M=2$ has an $\alpha=1.3390$ (as given by Formula~(7.1.23), Page~225 in \cite{daubechies:1992}\footnote{The $\alpha$ in the
table on Page~226 is our $\alpha$ minus 1.}), the condition on $\alpha$ is satisfied only for $H>0.67$.  
How should one proceed in a situation where $H\in(0,1)$ is unknown?
There are three alternatives: 1) Use  $M\geq3$ since the condition is satisfied for the Daubechies wavelet with $M\geq3$ (for
which $\alpha>1.63$). 2) Sum the FGN to get a discrete-time version of
FBM as above and use $M\geq2$. 3) Use $M=2$ and apply Theorem~\ref{theo:VjapproxDensityJbound} with a smoothness index
$\beta'<\beta$ instead of $\beta$, worsening the bound in~(\ref{eq:DjApprox}). 
\end{remark}

\begin{remark}
When $d<0$, it is not $M$ but $\alpha$ which should influence the choice of the wavelet. The more negative the value of $d$,
the higher the required value of $\alpha$. Recall that a high value of $\alpha$ corresponds to a fast decrease of $\hat{\psi}(\xi)$ as
$|\xi|\to\infty$. 
\end{remark}
%
%
\section{Analysis of the memory parameter estimator based on the regression of the wavelet variance}
\label{sec:abry-veitch}
In this section, we consider a \emph{Gaussian} process $X$ with memory parameter $d$ and generalized spectral density 
$f(\lambda)=|1- \rme^{-\rmi \lambda}|^{-2d} f^\ast(\lambda)$. Then, for any $K>(d-1/2)$, the distribution of the
$K$-th order increment process $\diffop^K X$ only depends on $d$ and $f^\ast$. 
We apply Theorem~\ref{theo:VjapproxDensityJbound} to
study the wavelet estimator of the memory parameter $d$, based on the regression of the scale spectrum
$\vj{j}{d,f^\ast}$ with respect to the scale index $j$.
This is reasonable because, for large scale $j$, $\log\vj{j}{d,f^\ast}$ is  approximately an affine function of
$j$ with slope $(2 \log2)\,d$ (see~(\ref{eq:VjApprox}) in Theorem \ref{theo:VjapproxDensityJbound}).
Given $n$ observations $X_1, \dots, X_n$, $\vj{j}{d,f^\ast}$
can be estimated by the empirical variance
\begin{equation}\label{eq:hvj}
\hvj{j}{n_j} \eqdef n_j^{-1} \sum_{k=0}^{n_j-1} \left( \dwt_{j,k}^X \right)^2 \eqsp ,
\end{equation}
where for any $j$, $n_j$ denotes the number of available wavelet coefficients at scale index $j$, namely, from~(\ref{eq:deltan}),
\begin{equation}\label{eq:ni}
n_j = [2^{-j} (n - \L + 1) - \L +1] \; ,
\end{equation}
where $n$ is the size of the time series and $[x]$ denotes the integer part of $x$.
An estimator of the memory parameter $d$ is then obtained by
regressing the logarithm of the empirical variance $\log (\hvj{i}{n_i})$  for a finite number of scale indices
$j \in \{J_0, \dots, J_0+\ell\}$
where $J_0$ is the lower scale and $1+\ell\geq2$ is the number of scales in the regression.
For a sample size equal to $n$, this estimator is well defined for $J_0$ and $\ell$ such that $\ell\geq1$
and $J_0+\ell\leq J(n)$ where
\begin{equation}
\label{eq:def:Jn}
J(n) \eqdef [\log_2 (n-\L+1)-\log_2(\L)]
\end{equation}
is the maximal index $j$ such that $n_j\geq1$. The regression estimator can be expressed formally as
\begin{equation}
\label{eq:definition:estimator:regression}
\hat{d}_n(J_0,\regressweights) \eqdef \sum_{j=J_0}^{J_0+\ell} w_{j-J_0} \log \left( \hvj{j}{n_j} \right) \eqsp ,
\end{equation}
where the vector $\regressweights \eqdef[w_0,\dots,w_{\ell}]^T$  satisfies
\begin{equation}
\label{eq:propertyw}
\sum_{i=0}^{\ell} w_{i}  = 0\quad\text{and}\quad 2 \log(2) \sum_{i=0}^{\ell} i w_{i}  = 1 \eqsp.
\end{equation}
One may choose, for example, $\regressweights$ corresponding to the weighted least-squares regression vector, defined by
\[
\regressweights= D B (B^TDB)^{-1} \mathbf{b}  \eqsp,
\]
where
$B \eqdef \left[\begin{matrix}
1 & 1 & \dots & 1 \\
0 & 1 & \dots & \ell
\end{matrix}\right]^T$ is the so-called design matrix, $D$ is a definite positive matrix and
\begin{equation}\label{eq:bfbDef}
\mathbf{b}\eqdef [0 \;\; (2\log(2))^{-1}]^T.
\end{equation}
Ordinary least square regression corresponds to the case where $D$ is the identity matrix.


We now compute a bound of the mean square error and an asymptotic equivalent of the variance of
$\hat{d}_n(J_0,\regressweights)$ in the usual semi-parametric framework adopted by Robinson and his co-authors for studying
Fourier estimators. For the wavelet estimator defined above, these quantities depend primarily on $n$ and on the scale index
$J_0$, while in the Fourier case, the bounds are generally expressed as functions of $n$ and a \emph{bandwidth} parameter
$m$, equal to the number of discrete Fourier frequencies used.
To ease comparison, we will express our results with respect to $n$ and $m$, where $m$ is the number of wavelet
coefficients appearing in $\hat{d}_n(J_0,\regressweights)$, namely,
$$
m \eqdef \sum_{j=J_0}^{J_0+\ell}n_j\eqsp.
$$
Since $\sum_{j=J_0}^{J_0+\ell}2^{-j}=2^{-J_0}(2-2^{-\ell})$, one gets immediately from~(\ref{eq:ni}) that
$|m-n2^{-J_0}(2-2^{-\ell})|\leq 2 (\ell+1) (\L-1)$. Thus 
$m\to\infty$ is equivalent to having $n2^{-J_0}\to\infty$, and, when these conditions hold, we have
\begin{equation}\label{eq:mSim}
m(n)\sim n2^{-J_0(n)}(2-2^{-\ell}) \eqsp .
\end{equation}
The next result provides a bound to the bias $\PE\left[ \hat{d}_n(J_0,\regressweights) \right] - d$
and to the variance $\PVar \left[ \hat{d}_n(J_0,\regressweights) \right]$.

\begin{theorem}
\label{theo:AbryVeitchEstimator}
Assume that \allWA\ hold with $M \geq 1$, $\alpha >1$ and that $X$ is Gaussian.
Let $\regressweights$ be a vector satisfying \eqref{eq:propertyw} for some $\ell \geq 1$.
Let  $\dmin$, $\dmax$ be two scalars such that $\dmin < \dmax$ and
$[\dmin,\dmax] \subset ( (1+\beta)/2-\alpha,M]$, where $\beta \in (0,2]$.
Then, there exists a finite constant $C$ (depending only
on $\regressweights$, $\beta$, $L$, $\dmin$, $\dmax$, $\phi$ and $\psi$) such that for all $J_0 \in \left\{0, \dots, J(n) -\ell \right\}$,
$d \in [\dmin, \dmax]$, and $f^\ast \in \calH(\beta,L)$ with $f^\ast(0)>0$
\begin{align}
\label{eq:bound:bias}
&\left| \PE \left[ \hat{d}_n(J_0,\regressweights) \right] - d \right| \leq C \left
\{\left(\frac{m}{n}\right)^{\beta} + m^{-1}\right\} \eqsp, \\
\label{eq:bound:variance}
&\PVar \left[ \hat{d}_n(J_0,\regressweights) \right] \leq C  \left\{ m^{-1} + \1\left(\frac {m}{n} \geq C^{-1}\right) \right\}\eqsp .
\end{align}
\end{theorem}

\begin{remark}
While the bias term bound contains $(m/n)^\beta$, the variance bound has an indicator function which is zero for sufficiently
small values of $m/n$, hence is $o((m/n)^{2\beta})$. This indicator function cannot be dispensed with. Indeed if we start our
estimation at the finest scale $J_0=0$, corresponding roughly to $m=n$, all we can say is that $\PVar(\hat{d}_n)\leq C$. If,
however, we start our estimation at a coarse enough scale $J_0$, corresponding to $m/n < C^{-1}$, then  
  $\PVar(\hat{d}_n)$ is bounded by $C m^{-1}$, which tends to zero as $m\to\infty$.
\end{remark}
By combining \eqref{eq:bound:bias} and \eqref{eq:bound:variance} it is possible to obtain a bound on the mean square error of
$\hat{d}_n (J_0,\regressweights)$.
More precisely, there exists a constant $C$ (depending only on $M$, $\alpha$, $\beta$, $L$, $\dmin$ and $\dmax$) such that,
for any $f^\ast \in \calH(\beta,L)$, $d \in [\dmin,\dmax]$ and $J_0 \in \left\{0, \dots, J(n)-\ell  \right\}$,
\begin{equation}\label{eq:MSE}
\PE \left[ \left\{ \hat{d}_n (J_0,\regressweights) -d \right\}^2 \right]\leq
C \left\{ \left(\frac{m}{n}\right)^{2\beta} + m^{-1} \right\} \eqsp.
\end{equation}
This shows in particular that, for any non-decreasing sequence $\{ J_0(n), n \geq 0\}$ such that
$ m^{-1}+ m/n\to0$, $\hat{d}_n(\regressweights) \eqdef \hat{d}_n(J_0(n),\regressweights)$ is a consistent estimator of $d$.
If the regularity exponent $\beta$ is known, it is possible to choose  $J_0(n)$
to balance these two terms, that is, set $(m/n)^{2\beta}\asymp m^{-1}$ or equivalently $2^{J_0(n)} \asymp n^{1/(1+2\beta)}$
as $n \to \infty$. If we choose $J_0(n)$ in such a
way,~(\ref{eq:mSim}) and~(\ref{eq:MSE}) imply
\[
\limsup_{n \to \infty}  \sup_{d \in [\dmin, \dmax]} \sup_{f^\ast \in \calH(\beta,L)} n^{2\beta/(1+2\beta)} \PE \left[ \left\{ \hat{d}_n (\regressweights) -d \right\}^2 \right]  < \infty \eqsp.
\]
As shown in \cite{giraitis:robinson:samarov:1997}, $n^{-2 \beta/(1+2 \beta)}$ is the
minimax rate of convergence for the memory parameter $d$ in this semi-parametric setting. Therefore,
\begin{corollary}
The wavelet estimator is rate optimal in the minimax sense.
\end{corollary}

We shall now obtain the asymptotic behavior of $\PVar \left[\hat{d}_n(\regressweights) \right]$ as $n\to\infty$.

\begin{theorem}
\label{theo:asympvariance}
Assume that \allWA\ hold with $M \geq 1$, $\alpha >1$  and that $X$ is Gaussian.
Let $\regressweights$ be a vector satisfying \eqref{eq:propertyw} for some $\ell \geq 1$.
Let $\{ J_{0}(n), n \in \Nset\}$ be a sequence such that  
$m\to \infty$ as $n \to \infty$.
For any $f^\ast \in \calH(\beta,L)$, where $\beta \in (0,2]$, and $ d \in ( (1+\beta)/2-\alpha,M]$,
\begin{equation}
\label{eq:limitingexpressionvariance}
\lim_{n \to \infty} m \PVar \left[\hat{d}_n(\regressweights) \right] =
(2-2^{-\ell}) \, \regressweights^T \AVvar[\psi]{d}{} \regressweights \eqsp,
\end{equation}
where $\AVvar[\psi]{d}{}$ is the $(1+\ell)\times(1+\ell)$ matrix defined as
\begin{align}
\label{eq:definitionAVvar}
\AVvar[\psi]{d}{i,j} \eqdef \frac{4 \pi 2^{2d |i-j| } 2^{i \wedge j}}{\Kvar[\psi]{d}^2} \int_{-\pi}^{\pi} \left|
  \bdensasymp[\psi]{|i-j|}{\lambda}{d} \right|^2 \, d \lambda  && 0 \leq i, j \leq \ell \eqsp .
\end{align}
\end{theorem}

\begin{remark}
The asymptotic expression of the variance \eqref{eq:limitingexpressionvariance} is a quadratic form of $\regressweights$
defined by the matrix $\AVvar[\psi]{d}{}$, which depends only on $d$ and $\psi$ (see \eqref{eq:definitionAVvar}).
The standard theory of linear regression shows that, for any $\ell\geq1$, the optimal regression vector of length $\ell+1$ is
\[
\regressweights^{\mathrm{opt}}(d,\psi) \eqdef  \AVvarInv[\psi]{d}{} B (B^T \AVvarInv[\psi]{d}{} B)^{-1}  \mathbf{b}
\]
and the associated limiting variance is $(2-2^{-\ell})\mathbf{b}^T(B^T \AVvarInv[\psi]{d}{} B)^{-1} \mathbf{b}$.
This optimal regression vector cannot be used directly since it depends on $d$ which is unknown, but one may apply a two-step procedure
using a preliminary estimate of $d$ as in \cite{bardet:2002} in a similar context.
\end{remark}

\begin{remark}
When computing confidence intervals in practice, one sometimes uses asymptotic variances
in~(\ref{eq:limitingexpressionvariance}) with $\AVvarInv[\psi]{0}{}$ instead of $\AVvarInv[\psi]{d}{}$, see
e.g.~\cite{abry:veitch:1998}. The expression $\AVvarInv[\psi]{0}{}$ can be easily obtained if the wavelets are
orthonormal. In this case, by~(\ref{eq:bdensasympNullD}) and~(\ref{eq:definitionAVvar}), for $i\neq j$, $\AVvar{0}{i,j}=0$
and $\AVvar{0}{j,j}=8\pi^2\,2^{j}/\Kvar{0}^2=2^{j+1}$ since by~(\ref{eq:Kpsi}), 
$\Kvar{0}=\int_\Rset|\hat{\psi}(\xi)|^2\,d\xi=2\pi\int_\Rset|\psi(t)|^2\,dt=2\pi$. 
Then~(\ref{eq:limitingexpressionvariance}) becomes 
$$
\lim_{n \to \infty} m \PVar \left[\hat{d}_n(\regressweights) \right] =
(2-2^{-\ell}) \, 2\, \sum_{j=0}^\ell w_j^2\,2^j \; .
$$
One can reformulate this in terms of $n_{J_0}\sim n2^{-J_0}$ instead of $m$. In view of~(\ref{eq:ni}) and~(\ref{eq:mSim}),
one gets the following simple expression of the asymptotic variance when $d=0$:
$$
\lim_{n \to \infty} n_{J_0(n)} \PVar \left[\hat{d}_n(\regressweights) \right] = 2\, \sum_{j=0}^\ell w_j^2\,2^j \; .
$$
\end{remark}

If we choose $m(n)$ (or $J_0(n)$) such that the bias in~(\ref{eq:bound:bias}) is asymptotically negligible, then we can obtain
the asymptotic behavior of the mean square error $ \PE \left( \hat{d}_n(\regressweights) - d \right)^2$. In view
of~(\ref{eq:bound:bias}) and~(\ref{eq:limitingexpressionvariance}), we need $m\to\infty$ and
$\{(m/n)^\beta+m^{-1}\}^2<<m^{-1}$, or equivalently
\begin{equation}
\label{eq:AsympCondOnmForCLT}
n 2^{-J_0(n)(1+2\beta)} + n^{-1} 2^{J_0(n)} \to 0,\quad n\to\infty \eqsp .
\end{equation}

\begin{corollary}
If~(\ref{eq:AsympCondOnmForCLT}) holds, then
for $f^\ast \in \calH(\beta,L)$ and $ d \in ( (1+\beta)/2-\alpha,M]$,
\[
\lim_{n \to \infty} n 2^{-J_0(n)} \PE \left( \hat{d}_n(\regressweights) - d \right)^2 =
\regressweights^T \AVvar[\psi]{d}{} \regressweights \eqsp.
\]
\end{corollary}
This result of course hints at the existence of a central limit theorem for the estimator $\hat{d}_n(\regressweights)$. Such
a result can be obtained by using a central limit theorem for quadratic forms of Gaussian variables which is established in a companion
paper Moulines, Roueff and Taqqu~(\citeyear{moulines:roueff:taqqu:2005b}).

\section{Proof of  Theorem \ref{theo:VjapproxDensityJbound}}
\label{sec:proof:prop:VjapproxDensityJbound}
From now on, we denote by $C$ constants possibly depending on $\dj$, $\dmin$, $\dmax$, $\beta$, $\phi$ and $\psi$, which may
change from line to line and we omit the dependence on $\phi$ and $\psi$ in the notations.
We assume, without loss of generality that $f^\ast(0)=1$.

\medskip

\noindent\textbf{Proof of} \ref{item:VjApprox}.
In the expression~(\ref{eq:definitionvj}) of $\vj{j}{d,f^{\ast}}$, $j\geq0$, we will approximate $|H_j(\lambda)|^2$
using~(\ref{eq:fj2Approx}). Thus define
$$
A_j\eqdef2^j\int_{-\pi}^\pi |1 - \rme^{-\rmi \lambda}|^{-2d} \, f^\ast(\lambda) \, |\hat{\phi}(\lambda)\hat{\psi}(2^j\lambda)|^2\,d\lambda
\quad\text{and}\quad
R_j \eqdef  \vj{j}{d,f^\ast} - A_j \eqsp .
$$
By~(\ref{eq:fj2Approx}), we have
\begin{equation}\label{eq:Rj}
|R_j|\leq C\,2^{j(1+M-\alpha)}\, \int_{-\pi}^\pi |1 - \rme^{-\rmi \lambda}|^{-2d} \, f^\ast(\lambda) \,|\lambda|^{2M}\,(1+2^j|\lambda|)^{-\alpha-M}\,d\lambda \eqsp.
\end{equation}
We consider $A_j$ and $R_j$ separately starting with $A_j$. 

Express $A_j$ as
\begin{equation}\label{eq:Aj}
A_j=2^j\int_{-\pi}^\pi g(\lambda) |\lambda|^{-2d}\,|\hat{\psi}(2^j\lambda)|^2\,d\lambda \eqsp,
\end{equation}
where
\begin{equation}\label{eq:g}
 g(\lambda) \eqdef f^\ast(\lambda)|\hat{\phi}(\lambda)|^2\,\left|\lambda/(1-\rme^{\rmi\lambda})\right|^{2d},
\quad \lambda\in(-\pi,\pi) \eqsp.
\end{equation}
Since $\hat{\phi}$ is infinitely differentiable by \ref{item:Wreg},
$\lambda \mapsto |\hat{\phi}(\lambda)|^2\,\left|\lambda/(1-\rme^{\rmi\lambda})\right|^{2d}$ is  infinitely differentiable on $[-\pi,\pi]$.
Because $f^\ast \in \calH(\beta,L)$ and $f^\ast(0)=1$, there exists a constant $C$
such that for all $\lambda\in(-\pi,\pi)$,
\begin{equation}\label{eq:gHolder}
|g(\lambda)-g(0)|\leq C\, L \, |\lambda|^\beta \eqsp ,
\end{equation}
where $g(0)=1$ because $\hat{\phi}(0)=1$ by Condition~\ref{item:Wreg}.
Moreover, $A_j$ is finite  by~\eqref{eq:MVM} since $g$ is bounded and $M>d-1/2$.
We shall now replace the function $g(\lambda)$  by the constant $g(0)= 1$ and
extend the interval of integration from $[-\pi,\pi]$ to the whole real line in~(\ref{eq:Aj}).
Eqs.~\eqref{eq:Aj} and~\eqref{eq:gHolder} imply
$$
\left|A_j-2^j\int_{-\pi}^\pi g(0)\,|\lambda|^{-2d}\,|\hat{\psi}(2^j\lambda)|^2\,d\lambda \right|
\leq C\, L \, 2^j \int_{-\pi}^\pi |\lambda|^{\beta-2d}\,|\hat{\psi}(2^j\lambda)|^2\,d\lambda \eqsp.
$$
First observe that, after a change of variable,
$$
2^j\int_{-\pi}^\pi |\lambda|^{\beta-2d}\,|\hat{\psi}(2^j\lambda)|^2\,d\lambda
\leq 2^{j(2d-\beta)} \int_{-\infty}^\infty \left\{ |\lambda|^{\beta-2\dmin} \vee |\lambda|^{\beta-2\dmax} \right\} \,|\hat{\psi}(\lambda)|^2\,d\lambda
$$
In the RHS of this inequality, using the behavior of $|\hat{\psi}(\lambda)|$ at infinity and at the origin implied by
\ref{item:psiHat} and ~\ref{item:MVM} respectively, and because $\dmax < M+1/2$ and $\dmin >(1+\beta)/2 - \alpha$, the integral is a
finite constant.
We further observe that, by \ref{item:psiHat}, since $\dmin>1/2-\alpha$, we may write
$$
2^j \int_{|\lambda|>\pi}|\lambda|^{-2d}\,|\hat{\psi}(2^j\lambda)|^2\,d\lambda
\leq C\, 2^{j(1-2\alpha)} \int_{|\lambda|>\pi} |\lambda|^{-2(\alpha+\dmin)}\,d\lambda 
$$
which is integrable.
Since, by~(\ref{eq:CondDminDmax-a}), $1-2\alpha<2d-\beta$, there exists a constant $C$ such that
\begin{equation}
\left|A_j-  \, \Kvar{d} \, 2^{2jd} \right| \leq \\
C\,L \, 2^{(2d-\beta)j} \eqsp,
\label{eq:Ajapprox}
\end{equation}
where $\Kvar{d} $ is given by~(\ref{eq:Kpsi}).

We now compute a bound for $R_j$ using~(\ref{eq:Rj}).
Note that
there exists a constant $C$
such that, for all $\lambda \in (-\pi,\pi)$,
\begin{equation}\label{eq:fBound}
f(\lambda)= f^\ast(\lambda)\,\left|\frac{\lambda}{1-\rme^{\rmi\lambda}}\right|^{2d}\,|\lambda|^{-2d}\leq
C\,L |\lambda|^{-2d}.
\end{equation}
Plugging this into \eqref{eq:Rj} and then separating $\lambda<1$ and $\lambda\geq1$, we obtain
\begin{align*}
R_j  &\leq C L 2^{2jd} 2^{-j(M+\alpha)} \int_{0}^{2^j \pi} \left\{ \lambda^{2(M-\dmin)} \vee \lambda^{2(M-\dmax)} \right\}  \left(1 + \lambda \right)^{-\alpha -M} d \lambda \\
    &\leq C L 2^{j(2d-\beta)} 2^{-j(M+\alpha-\beta)} \left\{ \int_0^1 \lambda^{2(M-\dmax)} d \lambda + \int_{1}^{2^j \pi} \lambda^{M-2\dmin-\alpha} d \lambda \right\}.
\end{align*}
Since $2(M-\dmax)>-1$, the first integral is a finite constant. Depending on whether $M-2\dmin-\alpha$ is less than, equal to
or larger than $-1$ the second integral is bounded by a finite constant, $\log\pi+j\log2$ or $C2^{j(1+M-2 \dmin-\alpha)}$. In
the two first cases, we simply observe that $M\geq1$, $\alpha>1$ and $\beta\leq2$ imply 
$M+\alpha-\beta > 0$, and in the last case that $-(M+\alpha-\beta)+1+M-2
\dmin-\alpha=1-2\dmin-2\alpha+\beta\leq0$ by~(\ref{eq:CondDminDmax-a})
so that, in all cases, $R_j\leq C\,L \, 2^{(2d-\beta)j}$. This condition, with (\ref{eq:Ajapprox}), implies 
$$
|\vj{j}{d,f^{\ast}}-\Kvar{d} 2^{2jd}|=|A_j+R_j-\Kvar{d} 2^{2jd}|\leq C\,L\,2^{(2d-\beta)j} 
$$
which proves~\eqref{eq:VjApprox}.

\medskip

\noindent\textbf{Proof of} \ref{item:DjApprox}.
For ease of notation, we only consider the case $\dj=0$ so that $\be_\dj(\xi)=1$.
It is also enough to suppose $j \geq 1$. In \eqref{eq:definitionfj}, the summands are $2^j (2\pi)$-periodic; hence, omitting the summands,
$\sum_{l=0}^{2^j-1} = \sum_{l=0}^{2^{j-1}-1} + \sum_{j=2^{j-1}}^{2^j-1} = \sum_{l=0}^{2^{j-1}-1} + \sum_{l=-2^{j-1}}^{-1} =
\sum_{l=-2^{j-1}}^{2^{j-1}-1}$. Note that,  for $l \in \{-2^{j-1}, \dots, 2^{j-1}-1\}$ and $\lambda \in (0,\pi)$, we have 
$2^{-j}(\lambda+2 l \pi) \in (-\pi,\pi)$  so that  \eqref{eq:fj2Approx} applies.
Hence, $\bdens{j,0}{\lambda}{d,f^\ast}$ in~(\ref{eq:definitionfj}) is expressed as the sum of two functions
$A_j(\lambda)+R_j(\lambda)$, defined for all $\lambda\in(0,\pi)$ by
\begin{equation}\label{eq:Alambdaj}
A_j(\lambda)\eqdef
\sum_{l=-2^{j-1}}^{2^{j-1}-1} |2^{-j}(\lambda+2l\pi)|^{-2d}\,g(2^{-j}(\lambda+2l\pi))\,
|\hat{\psi}(\lambda+2l\pi)|^2
\end{equation}
where $g$ is defined in \eqref{eq:g} and where by~(\ref{eq:fj2Approx}),
\begin{equation}\label{eq:Rlambdaj}
R_j(\lambda)\leq C\,L\,2^{j(2d-M-\alpha)}
\sum_{l=-2^{j-1}}^{2^{j-1}} |\lambda+2l\pi|^{2(M-d)}\,(1+|\lambda+2l\pi|)^{-\alpha-M}.
\end{equation}
From~(\ref{eq:gHolder}), we get, for all $\lambda\in(0,\pi)$,
\begin{equation}\label{eq:Bl}
\left|A_j(\lambda)-2^{2dj}\,g(0)\,
\sum_{l=-2^{j-1}}^{2^{j-1}-1} |\lambda+2l\pi|^{-2d}\,|\hat{\psi}(\lambda+2l\pi)|^2\right| 
\leq C\, L \,2^{(2d-\beta)j}\,B_j(\lambda),
\end{equation}
where, by~(\ref{eq:MVM}) and~\ref{item:psiHat}, for all $\lambda\in(0,\pi)$,
\begin{align}
\nonumber
B_j(\lambda)&\eqdef\sum_{l=-2^{j-1}}^{2^{j-1}-1} |\lambda+2l\pi|^{\beta-2d}\,|\hat{\psi}(\lambda+2l\pi)|^2\\
\nonumber
&\leq C\,\left(|\lambda|^{\beta+2(M-d)} +2\sum_{l\geq1}|\lambda+2l\pi|^{\beta-2d-2\alpha}\right)\\
\label{eq:BjBound}
&\leq C\,\left(1+2\sum_{l\geq1}(2l-1)^{\beta-2\dmin-2\alpha}\right) <\infty
\end{align}
since $|\lambda+2l\pi|\geq \pi(2l-1)$, $\beta>0$ and, by~(\ref{eq:CondDminDmax-b}), one has $M\geq d$ and $\beta-2 \dmin
-2\alpha<-1$. By the same arguments, for all $\lambda\in(0,\pi)$,
$$
\sum_{|l|\geq 2^{j-1}-1} |\lambda+2l\pi|^{-2d}\,|\hat{\psi}(\lambda+2l\pi)|^2
\leq  C 2^{j(1-2(\dmin+\alpha))} 
$$
is bounded since the exponent is negative.
Eqs.~(\ref{eq:bDpsi}) with $\dj=0$, (\ref{eq:Bl}), $g(0)=1$ and the above inequalities yield that, for all
$\lambda\in(0,\pi)$,
$$
\left|A_j(\lambda)-\,\bdensasymp{0}{\lambda}{d}\,2^{2dj}\right|
\leq C\,L\, 2^{(2d-\beta)j}\, \eqsp.
$$

We now turn to bounding $R_j(\lambda)$ using~(\ref{eq:Rlambdaj}). For all $\lambda\in(0,\pi)$, using $|\lambda-2
l\pi|\geq\pi(2l-1)$ and~(\ref{eq:CondDminDmax-b}),
\[
R_j(\lambda) \leq C\, L\,2^{j(2d-\beta)} 2^{-j(M+\alpha-\beta)}\, \left(1+\sum_{l=1}^{2^{j}} l^{-2 \dmin+M-\alpha}\right)
\]
which can be bounded as in the proof of~\ref{item:VjApprox}, by considering the cases 
$M-2\dmin-\alpha<$, $=$ or $>-1$.

The joint continuity of $(\lambda,d)\mapsto\bdensasymp[\psi]{0}{\lambda}{d}$ on $\Rset\times [\dmin,\dmax]$  follows
from~(\ref{eq:BjBound}) and dominated convergence.

\section{Proofs of Theorem~\ref{theo:AbryVeitchEstimator} and  Theorem~\ref{theo:asympvariance}}\label{sec:proofs-theor-refth}
From now on, we denote by $J_{\min}$,
$n_{\min}$, $C$ and $C'$ some positive constants whose values may change upon each occurrence and which depend at most on 
$\regressweights$, $\beta$, $L$, $\dmin$, $\dmax$, $\phi$, and $\psi$. We will repeatedly use that,
by~(\ref{eq:ni}),~(\ref{eq:def:Jn}) and~(\ref{eq:mSim}), for $J_0(n)\in\{0,\dots,J(n)-\ell\}$,  
\begin{equation}\label{eq:asympmetc} 
n_{J_0(n)}\asymp n_{J_0(n)+\ell}\asymp n2^{-J_0(n)}\asymp m(n),
\end{equation}
where $a \asymp b$ means that there is a constant $C\geq1$ such that $a/C\leq b \leq C\,a$. 
Finally, for any measurable vector-valued function $\varphi$ on $[-\pi,+\pi]$ and any $p > 0$, 
$\| \varphi \|_p = \left( \int_{-\pi}^{\pi} |\varphi(\lambda)|^p d \lambda  \right)^{1/p}$.
\begin{proposition}
\label{prop:varlimhvj}
Let $\hvj{j}{n_j}$ be defined as in~(\ref{eq:hvj}) and $\bdensasymp[\psi]{\dj}{\lambda}{d}$, $\dj\geq0$, be defined as
in~(\ref{eq:bDpsi}). Under the assumptions of 
Theorem~\ref{theo:AbryVeitchEstimator}, one has, as $j\to\infty$ and $n_j\to\infty$, 
$$
2^{-4dj}\,n_{j-\dj}\,\PCov\left(\hvj{j}{n_j},\hvj{j-\dj}{n_{j-\dj}}\right) \to (f^\ast(0))^2\, 4\pi \,
\|\bdensasymp[\psi]{\dj}{\cdot}{d}\|_2^2 \;,  
$$
uniformly on $d\in[\dmin,\dmax]$ and $f^\ast\in\calH(\beta,L)$.
\end{proposition}
\begin{proof}
We set $f^\ast(0)=1$ without loss of generality. Using~(\ref{eq:hvj}) and~(\ref{eq:Defbd}), we write
\begin{align}
\nonumber
\PCov
\left[ \hvj{j}{n_j}, \hvj{j-\dj}{n_j} \right]  &= \frac{1}{n_{j} n_{j-\dj}} \sum_{k,l=0}^{n_{j}-1} \sum_{\dk=0}^{2^{\dj}-1}
\PCov\left[\dwt^2_{j,k}, \dwt^2_{j-\dj,l 2^{\dj}+\dk} \right] \\
\label{eq:borne:covariance:warianceempirique}
& = \frac{2}{n_{j} n_{j-\dj}} \sum_{k,l=0}^{n_{j}-1} \left| \PCov\left[ \dwt_{j,k}, \bdwt_{j,l}(\dj) \right] \right|^2 \\
\label{eq:borne:covariance:warianceempiriqueTau}
& = \frac{2}{n_{j-\dj}} \sum_{\tau\in\Zset} \left(1-\frac{|\tau|}{n_{j}}\right)_+ \left| \PCov\left[ \dwt_{j,0}, \bdwt_{j,\tau}(\dj) \right]
\right|^2 \eqsp,
\end{align}
where, in~(\ref{eq:borne:covariance:warianceempirique}), we used the fact that if the scalar $X$ and the vector
$\bY=[Y_1\,\,\dots\,\,Y_p]^T$ are jointly Gaussian, 
$$
\sum_{k=1}^p\PCov\left(X^2,Y_k^2\right)=2\sum_{k=1}^p\PCov^2(X,Y_k)=2\left|\PCov(X,\bY)\right|^2.
$$
Using the notation $M_n$ defined in~(\ref{eq:Nn}), we have
\begin{equation}\label{eq:NnBdens}
\left(\sum_{\tau\in\Zset} \left(1-\frac{|\tau|}{n_{j}}\right)_+ \left| \PCov\left( \dwt_{j,0}, \bdwt_{j,\tau}(\dj) \right) \right|^2\right)^{1/2} = 
M_{n_j}(\bdens{j,\dj}{\cdot}{d,f^\ast}) \; ,
\end{equation}
since, by Corollary~\ref{cor:bD}, $\bdens{j,\dj}{\cdot}{d,f^\ast}$ is the cross-spectral density of the vector
$[\dwt_{j,0},\,\bdwt_{j,\tau}(\dj)]$. 
Applying Lemma~\ref{lem:empVarL2bound}--(\ref{eq:NnLipschitz}), the relation $\|\cdot\|_2\leq \sqrt{2\pi}\|\cdot\|_\infty$ and
Theorem~\ref{theo:VjapproxDensityJbound}--(\ref{eq:DjApprox}), there is a constant $C$ such that
\begin{equation}\label{eq:NnBdensAsymp}
\left| M_{n_j}(\bdens{j,\dj}{\cdot}{d,f^\ast})
-2^{2jd}\,M_{n_j}(\bdensasymp{\dj}{\cdot}{d})\right| \leq C \, 2^{(2d-\beta)j} \; .
\end{equation}
On the other hand, by Lemma~\ref{lem:empVarL2bound}--(\ref{eq:NnLim}), we have, as $n_j\to\infty$, 
\begin{equation}\label{eq:NnBdensAsympLim}
|M_{n_j}(\bdensasymp{\dj}{\cdot}{d})|^2 \to 2\pi\,\|\bdensasymp{\dj}{\cdot}{d}\|_2^2 \; .
\end{equation}
The convergence in~(\ref{eq:NnBdensAsympLim}) holds uniformly on $d\in[\dmin,\dmax]$ because of the joint continuity of
$(\lambda,d)\mapsto\bdensasymp{\dj}{\lambda}{d}$ stated in Theorem~\ref{theo:VjapproxDensityJbound}.

The result follows from~(\ref{eq:borne:covariance:warianceempiriqueTau})--(\ref{eq:NnBdensAsympLim}).
\end{proof}

\begin{proof}[Proof of Theorem \ref{theo:AbryVeitchEstimator}]
Again, we set $f^\ast(0)=1$ without loss of generality. 
The bias $\PE [ \hat{d}_n(J_0,\regressweights)]-d$ can be decomposed into two terms as follows
\begin{multline}
\label{eq:decomposition-biais}
\sum_{j=J_0}^{J_0+\ell} w_{j-J_0} \PE \left[ \log ( \hvj{j}{n_j}) \right] - d = \sum_{j=J_0}^{J_0+\ell} w_{j-J_0} \log \left[
  \vj{j}{d,f^\ast} \right]  - d + \\
\sum_{j=J_0}^{J_0+\ell} w_{j-J_0} \left\{ \PE \left[ \log ( \hvj{j}{n_j})  \right] - \log \left[ \PE[\hvj{j}{n_j}] \right] \right\} \eqsp,
\end{multline}
where $\hvj{j}{n_j}$ is the wavelet coefficient empirical variance~(\ref{eq:hvj}) and $\PE[\hvj{j}{n_j}]=\vj{j}{d,f^\ast}$.

Using \eqref{eq:propertyw}, the first term on the RHS of~(\ref{eq:decomposition-biais}) may be rewritten as
\begin{equation}
\label{eq:bais01}
\sum_{j=J_0}^{J_0+\ell} w_{j-J_0} \log \left[ \vj{j}{d,f^\ast} \right]  - d 
= \sum_{j=J_0}^{J_0+\ell} w_{j-J_0} \log \left( 1 + \frac{\vj{j}{d,f^\ast} - \Kvar{d}\, 2^{2jd}}
{ \Kvar{d} \, 2^{2jd}} \right)
\end{equation}
By Theorem \ref{theo:VjapproxDensityJbound}-\eqref{eq:VjApprox} and using that 
$\inf_{d\in[\dmin,\dmax]} \Kvar{d} >0$, there exists a constant $C$ such that
$$
\frac{|\vj{j}{d,f^\ast} - \Kvar{d}\, 2^{2jd}|}{\Kvar{d} \, 2^{2jd}} \leq C 2^{-\beta j}.
$$ 
Using that $|\log(1+x)|\leq2|x|$ for $x\in(-1/2,\infty)$, there is a $J_{\min}$ such that, for $j\geq J_{\min}$, the
logarithm in the RHS of~(\ref{eq:bais01}) is bounded by $C\,2^{-\beta j}$, and by~(\ref{eq:asympmetc}), for all $J_0\geq
J_{\min}$, 
\begin{equation}
\label{eq:bias1}
\left|\sum_{j=J_0}^{J_0+\ell} w_{j-J_0}  \log [ \vj{j}{d,f^\ast} ] - d\right|
\leq C\,\sum_{j=0}^{\ell} |w_{j}| 2^{- (J_0+j) \beta}
\leq C \, \left(\frac mn\right)^\beta \eqsp .
\end{equation}
This bound is in fact valid for all $J_0\geq0$ because $\vj{j}{d,f^\ast}$ is bounded away from zero and infinity
independently of $d$ and $f^\ast$. Indeed, by~(\ref{eq:definitionvj}) and since $f^\ast\in\calH(\beta,L)$ with $f^\ast(0)=1$,
there is a small enough 
$\epsilon>0$ only depending on $\dmin$, $\dmax$, $\beta$ and $L$ such that
\begin{multline}\label{eq:vjUniversalBound}
C\,\int_\epsilon^{2\epsilon}|H_j(\lambda)|^2\,d\lambda \leq 
\int_{-\pi}^\pi|1-\rme^{-\rmi\lambda}|^{-2d}(1-L|\lambda|^\beta)_+|H_j(\lambda)|^2\,d\lambda
\leq \vj{j}{d,f^\ast} \\
\leq \int_{-\pi}^\pi|1-\rme^{-\rmi\lambda}|^{-2d}(1+L|\lambda|^\beta)_+|H_j(\lambda)|^2\,d\lambda
\leq C' \, \int_{-\pi}^\pi|\lambda|^{-2\dmax}|H_j(\lambda)|^2\,d\lambda \; .
\end{multline}
Observe that the lower bound in the previous display does not vanish since, as stated in Remark~\ref{rem:j0}, $H_j(\lambda)$
is a non-zero trigonometric polynomial for all 
$j\geq0$ and that the upper bound is finite since, by~(\ref{eq:Htilde}), $H_j(\lambda)=O(|\lambda|^M)$. Hence there is a
positive constant $C$ such that, for all $j=\{0,\dots,J_{\min}\}$,   
$C^{-1} \leq \vj{j}{d,f^\ast} \leq C$.   

We now consider the second term in the RHS of the display \eqref{eq:decomposition-biais}. The empirical
variance~(\ref{eq:hvj}) is a quadratic form in the wavelet coefficients at
$[\dwt_{j,0}, \dots, \dwt_{j,n_j-1}]$. By  Corollary~\ref{cor:bD}, these have spectral density
$\bdens{j,0}{\cdot}{d,f^\ast}$, given in \eqref{eq:definitionfj}.
By Lemma~\ref{lem:CovRadBound}, the spectral radius of the covariance matrix $\Gamma_{j}(d,f^\ast)$
of the random vector $[\dwt_{j,0}, \dots, \dwt_{j,n_j-1}]$
is bounded by the supremum of the spectral density,
\begin{equation}
\label{eq:bound:covariance:within-scale}
\rho \left[ \Gamma_{j}(d,f^\ast) \right]\leq 2 \pi \, \|\bdens{j,0}{\cdot}{d,f^\ast}\|_\infty \eqsp.
\end{equation}
Applying Proposition~\ref{prop:meanlogquadgauss}-\eqref{eq:boundmean} with $A=n_j^{-1}I_{n_j}$ and
$\Gamma=\Gamma_{j}(d,f^\ast)$ and using~(\ref{eq:bound:covariance:within-scale}), we get
\begin{equation}
\label{eq:bias2}
\left| \PE \left[ \log( \hvj{j}{n_j} ) \right] - \log \left[ \PE \left( \hvj{j}{n_j} \right) \ \right]\right| 
\leq 4\pi^2\,C\,\left(1\wedge n_j^{-1}\,\frac{\|\bdens{j,0}{\cdot}{d,f^\ast}\|_\infty^2}{n_j\,\PVar[\hvj{j}{n_j}]}\right)\eqsp,
\end{equation}
where $C$ is a universal constant. Now, by Theorem~\ref{theo:VjapproxDensityJbound}--(\ref{eq:DjApprox}) and by joint
continuity of $\bdensasymp{0}{\lambda}{d}$,  
$$
2^{-2dj}\|\bdens{j,0}{\cdot}{d,f^\ast}\|_\infty\leq C \, (1+\|\bdensasymp{0}{\cdot}{d}\|_\infty) \leq C \; .
$$
It follows from Proposition~\ref{prop:varlimhvj} that for $j\geq J_{\min}$ and $n_{j}\geq
n_{\min}$,
$$
2^{-4dj}\,n_j\,\PVar\left(\hvj{j}{n_j}\right) \geq 
2\pi\,\inf_{d\in[\dmin,\dmax]}\|\bdensasymp{0}{\cdot}{d}\|_2^2 
$$
which is positive by Remark~\ref{rem:BdAsympLpNorm}. The last two displayed equations imply that for $j\geq J_{\min}$ and $n_{j}\geq
n_{\min}$,
\begin{equation}
\label{eq:bound:LinftyOverVar}
\frac{\|\bdens{j,0}{\cdot}{d,f^\ast}\|_\infty^2}{n_j\,\PVar[\hvj{j}{n_j}]} \leq C  \; .
\end{equation}
Inserting~(\ref{eq:bound:LinftyOverVar}) into~(\ref{eq:bias2}) and using~(\ref{eq:asympmetc}), we get that for $J_0\geq
J_{\min}$ and $n_{J_0+\ell}\geq n_{\min}$, and $j=J_0,\dots,J_0+\ell$, 
\begin{equation}
\label{eq:biais2}
\left| \PE \left[ \log( \hvj{j}{n_j} ) \right] - \log \left( \PE \left[ \hvj{j}{n_j} \right] \ \right)\right| 
\leq C \, n_{j}^{-1}\leq C \, n_{J_0+\ell}^{-1} 
\leq C\, m^{-1} \leq C\,(m^{-1} + (m/n)^{\beta})\;.
\end{equation}
This last bound holds in fact without the preceding restrictions on $J_0$ and $n_{J_0+\ell}$. To see this,
use~(\ref{eq:bias2}) (with the ``bound 1'') and observe that, by~(\ref{eq:asympmetc}), $J_0\leq J_{\min}$ implies
$2^{-J_0}\geq2^{-J_{\min}}$, that is, $m/n\geq C$, and $n_{J_0+\ell}\leq n_{\min}$ implies $m^{-1}\geq C$. 
The bounds~\eqref{eq:bias1} and~(\ref{eq:biais2}), inserted in~(\ref{eq:decomposition-biais}), yield the bound~\eqref{eq:bound:bias} on
the bias. 

We now compute the variance of the estimator $\hat{d}_n(J_0,\regressweights)$. By
Proposition~\ref{prop:meanlogquadgauss}--(\ref{eq:boundcov}) and 
using \eqref{eq:bound:covariance:within-scale} and~(\ref{eq:bound:LinftyOverVar}) as in~(\ref{eq:bias2}), 
\begin{align}
\nonumber
&\left| \PVar( \hat{d}_n(J_0,\regressweights) ) - \sum_{i,j=J_0}^{J_0+\ell} w_{i-J_0} w_{j-J_0} 
\frac{\PCov(\hvj{i}{n},\hvj{j}{n})}
{\PE \left[ \hvj{i}{n_i} \right] \, \PE \left[ \hvj{j}{n_j} \right]}\right|  \\
\nonumber & \, \leq \sum_{i,j=J_0}^{J_0+\ell} | w_{i-J_0} w_{j-J_0} |
\left| \PCov\left( \log(\hvj{i}{n}) , \log(\hvj{j}{n})
\right) -
  \frac{\PCov(\hvj{i}{n},\hvj{j}{n})}{\PE \left[ \hvj{i}{n_i} \right]\PE \left[ \hvj{j}{n_j} \right]} \right| \\
\nonumber
& \, \leq C \sum_{i,j=J_0}^{J_0+\ell} | w_{i-J_0} w_{j-J_0} |
\left \{ 
\frac{n_i^{-3}\| \bdens{i,0}{\cdot}{d,f^\ast} \|_\infty^{3}}{\PVar^{3/2}(\hvj{i}{n})}
\vee  
\frac{n_j^{-3}\| \bdens{j,0}{\cdot}{d,f^\ast} \|_\infty^{3}}{\PVar^{3/2}(\hvj{j}{n})}
 \right\} \\
 \label{eq:VarsansLog}
 & \,  \leq C n_{J_0+\ell}^{-3/2} \leq C m^{-3/2} = o (m^{-1}) \eqsp.
\end{align}
On the other hand, by Proposition~\ref{prop:varlimhvj} and Theorem~\ref{theo:VjapproxDensityJbound}--(\ref{eq:VjApprox}), we
have, for any $\dj\geq0$, as $j_0\to\infty$ and $n_{j_0}\to\infty$,  
\begin{equation} 
\label{eq:LimCovVarEmp}
\frac{n_{j_0-\dj}\,\PCov(\hvj{j_0}{n_{j_0}},\hvj{j_0-\dj}{n_{j_0}})}
{\PE \left[ \hvj{j_0}{n_{j_0}} \right] \, \PE \left[ \hvj{j_0-\dj}{n_{j_0-\dj}} \right]} = 
\frac{n_{j_0-\dj}\,\PCov(\hvj{j_0}{n_{j_0}},\hvj{j_0-\dj}{n_{j_0}})}{\vj{j_0}{d,f^\ast} \vj{j_0-\dj}{d,f^\ast}} \to 
\frac{4\pi\,\|\bdensasymp{\dj}{\cdot}{d}\|_2^2}{2^{-2d\dj}\,(\Kvar{d})^2} 
\end{equation}
uniformly in $d\in[\dmin,\dmax]$ and $f^\ast\in\calH(\beta,L)$. Applying~(\ref{eq:LimCovVarEmp}) with $j_0=i\vee j$ and
$\dj=|i-j|$ so that $j_0-\dj=i\wedge j$, and since $\inf_{d\in[\dmin,\dmax]}\Kvar{d}>0$, for all $J_0\geq J_{\min}$ and
$n_{J_0+\ell}\geq n_{\min}$,  
$$
\left|\sum_{i,j=J_0}^{J_0+\ell} w_{i-J_0} w_{j-J_0}
\frac{\PCov[\hvj{i}{n},\hvj{j}{n}]}{\vj{i}{d,f^\ast} \vj{j}{d,f^\ast}} \right| \leq C \, n_{J_0+\ell}^{-1} \leq C \, m^{-1}\eqsp.
$$
This bound with \eqref{eq:VarsansLog} yields~\eqref{eq:bound:variance}  for $J_0\geq J_{\min}$ and $n_{J_0+\ell}\geq
n_{\min}$. When $J_0\leq J_{\min}$ or $n_{J_0+\ell}\leq n_{\min}$, the RHS of~\eqref{eq:bound:variance} is larger than a
positive constant  
and it suffices to use that, by the Minkowski inequality,
$$
 \PVar^{1/2}( \hat{d}_n(J_0,\regressweights) ) \leq 
\sum_{j=J_0}^{J_0+\ell} 
|w_{j-J_0}| \PVar^{1/2}\left( \log(\hvj{j}{n}) \right) \leq C \; ,
$$
where we applied Proposition~\ref{prop:meanlogquadgauss}--(\ref{eq:boundvar}) and that, by~(\ref{eq:vjUniversalBound}), 
$\PE[\hvj{j}{n}]=\vj{j}{d,f^\ast}$ does not vanish. 
\end{proof}

\begin{proof}[Proof of Theorem~\ref{theo:asympvariance}]
By~(\ref{eq:ni}) and~(\ref{eq:mSim}), for $i,j=J_0,\dots,J_0+\ell$,
$$
n_{i\wedge j}^{-1}\sim (n2^{- (i\wedge j)})^{-1}\sim m^{-1}(2-2^{-l})\,2^{(i-J_0)\wedge(j-J_0)}.
$$
Hence, by~(\ref{eq:VarsansLog}) and~(\ref{eq:LimCovVarEmp}),
$$
\lim_{m\to\infty} m \PVar( \hat{d}_n(J_0,\regressweights) ) = 
(2-2^\ell)
\sum_{i,j=J_0}^{J_0+\ell} w_{i-J_0} w_{j-J_0}
\frac{4\pi2^{2d|i-j|}}{(\Kvar{d})^2}\,2^{(i-J_0)\wedge(j-J_0)}\,\|\bdensasymp{|i-j|}{\cdot}{d}\|_2^2 \; ,
$$
which gives~(\ref{eq:limitingexpressionvariance}) after a change of variables.
\end{proof}

\appendix

\section{Approximation of wavelet filter transfer functions}\label{sec:appr-wavel-filt}

\begin{proposition}\label{prop:NewFjApprox}
Under \allWA, there exist positive constants $C_i$, $i=1,\dots,4$ only depending on $\phi$ and $\psi$,
such that, for all $j\geq0$ and $\lambda\in(-\pi,\pi)$,
\begin{eqnarray}
\label{eq:fjApprox}
|H_j( \lambda)-2^{j/2}\hat{\phi}(\lambda)\overline{\hat{\psi}(2^j\lambda)}|
& \leq & C_1\,2^{j(1/2-\alpha)}\,|\lambda|^M, \\
\label{eq:fjApproxBound}
|\hat{\phi}(\lambda){\hat{\psi}(2^j\lambda)}|
& \leq & C_2\,|2^j\lambda|^M\,(1+2^j|\lambda|)^{-\alpha-M},\\
\label{eq:fjBound}
|H_j( \lambda)|
& \leq & C_3\, 2^{j/2} \, |2^j\lambda|^M\,(1+2^j|\lambda|)^{-\alpha-M}, \\
\label{eq:fj2Approx}
\left|  |H_j( \lambda)|^2 - 2^{j}\,|\hat{\phi}(\lambda)\hat\psi(2^j\lambda)|^2\right|
& \leq & C_4\,2^{j(1+M-\alpha)}\,|\lambda|^{2M}\,(1+2^j|\lambda|)^{-\alpha-M}.
\end{eqnarray}
\end{proposition}
\begin{proof}
Under~\ref{item:Wreg} and~\ref{item:psiHat}, we have that, for all $t\in\Rset$,
$\sum_{k\in\Zset}\hat{\phi}(\lambda+2k\pi)\,\rme^{\rmi t (\lambda+2k \pi)}$ is a $2\pi$-periodic function, integrable on
$(-\pi,\pi)$ and whose $l$-th Fourier coefficients is
$$
\int_{-\pi}^\pi \sum_{k\in\Zset}\hat{\phi}(\lambda+2k\pi)\,\rme^{\rmi t(\lambda+2k \pi)}\,\rme^{-\rmi\lambda l}\,d\lambda
=\int_{-\infty}^\infty\hat{\phi}(\lambda)\,\rme^{\rmi t\lambda}\,\rme^{-\rmi \lambda l}\,d\lambda
=2\pi\,\phi(t-l).
$$
It follows that, for all $\lambda$ and $t$ in $\Rset$,
$$
\sum_{l\in\Zset}\phi(t-l)\,\rme^{\rmi \lambda l} = \sum_{k\in\Zset}\hat{\phi}(\lambda+2k\pi)\,\rme^{\rmi t (\lambda+2k \pi)},
$$
which is a form of the Poisson summation formula. Inserting this in~(\ref{eq:fjdef}) gives
\begin{eqnarray*}
H_j(\lambda)&=&2^{-j/2} \int_{-\infty}^\infty \left(\sum_{k\in\Zset}\hat{\phi}(\lambda+2k\pi)\,
\rme^{\rmi t(\lambda+2k \pi)}\right)\psi(2^{-j}t)\,dt\\
&=&2^{-j/2} \sum_{k\in\Zset}\hat{\phi}(\lambda+2k\pi)\int_{-\infty}^\infty \rme^{\rmi t (\lambda+2k\pi)}\psi(2^{-j}t)\,dt\\
&=&2^{j/2} \sum_{k\in\Zset}\hat{\phi}(\lambda+2k\pi)\overline{\hat{\psi}(2^j(\lambda+2k\pi))}.
\end{eqnarray*}
From this expression of $H_j$, we get, for all $j\geq0$ and
$\lambda\in(-\pi,\pi)$,
\begin{equation}\label{eq:fjeq}
|H_j( \lambda)-2^{j/2}\hat{\phi}(\lambda)\overline{\hat{\psi}(2^j\lambda)}|
= 2^{j/2} \left| \sum_{|k|\geq1}\hat{\phi}(\lambda+2k\pi)\overline{\hat{\psi}(2^j(\lambda+2k\pi))} \right|.
\end{equation}
Now using successively~(\ref{eq:MIM}) and~\ref{item:psiHat}, there is a constant $C$ such that, for all non-zero integer $k$
and all $\lambda\in(-\pi,\pi)$,
$ |\hat{\phi}(\lambda+2k\pi)|\leq C|\lambda|^M$ and
$$
|\overline{\hat{\psi}(2^j(\lambda+2k\pi))}|\leq C\,(2^j|\lambda+2k\pi|)^{-\alpha}
\leq C\,2^{-\alpha j}\,(2|k|\pi-|\lambda|)^{-\alpha}\leq \frac{C\,2^{-\alpha j}} {\pi^{\alpha}(2|k|-1)^{\alpha}}.
$$
Inserting these bounds into~(\ref{eq:fjeq}) gives~(\ref{eq:fjApprox}).

The bound~(\ref{eq:fjApproxBound}) follows from ~\ref{item:Wreg} ($|\hat\phi(\xi)| \leq \int_{-\infty}^\infty |\phi(t)| d t<\infty$),
~\ref{item:psiHat} ($|\hat{\psi}(\xi)| \leq C (1 + |\xi|)^{-\alpha}$) and~(\ref{eq:MVM}).

The two last bounds~(\ref{eq:fjBound}) and~(\ref{eq:fj2Approx}) follow from the two first~(\ref{eq:fjApprox})
and~(\ref{eq:fjApproxBound}).
Indeed, let $H^{(0)}_j(\lambda) \eqdef 2^{j/2}\,\hat{\phi}(\lambda)\overline{\hat\psi(2^j\lambda)}$.
For~(\ref{eq:fjBound}) we write
$$
|H_j( \lambda)|\leq |H_j(\lambda)-H^{(0)}_j( \lambda)| + |H^{(0)}_j( \lambda)|.
$$
Applying~(\ref{eq:fjApprox}) and~(\ref{eq:fjApproxBound}), the RHS of this equation is bounded by
\begin{multline*}
C_1\,2^{j(1/2-\alpha)}\,|\lambda|^M+C_2\,2^{j/2} |2^j\lambda|^M\,(1+2^j|\lambda|)^{-\alpha-M}
\leq \\
2^{j/2} |2^j\lambda|^M\,(1+2^j|\lambda|)^{-\alpha-M}
\,(C_1\,2^{-j(\alpha+M)}\,(1+2^j|\lambda|)^{\alpha+M}+C_2).
\end{multline*}
By observing that, for all $j\geq0$ and $\lambda\in(-\pi,\pi)$,
 the last term in parentheses is bounded by
$C_1\,2^{-j(\alpha+M)}\,(2^{1+j}\pi)^{\alpha+M}+C_2\leq C_1\,(2\pi)^{\alpha+M}+C_2$, we
get~(\ref{eq:fjBound}). For~(\ref{eq:fj2Approx}), we write
$$
\left|  |H_j( \lambda)|^2 - |H^{(0)}_j(\lambda)|^2 \right|  \leq
\left| H_j( \lambda) - H^{(0)}_j(\lambda)\right|
\left( | H^{(0)}_j( \lambda) | +  | H_j(\lambda)| \right)
$$
and apply~(\ref{eq:fjApprox}),~(\ref{eq:fjApproxBound}) and~(\ref{eq:fjBound}).
\end{proof}

\begin{corollary}
\label{cor:InfHj}
Under \allWA, there exists $j_0\geq0$ such that, for all $j\geq j_0$, $H_j$ is not identically zero.  
\end{corollary}
\begin{proof}
By~\ref{item:Wreg}, there exist sufficiently small positive constants $\epsilon$ and $\eta$ such that
$|\hat{\phi}(\lambda)|\geq1/2$ for all $|\lambda|\leq\epsilon$ and
$\inf_{|\lambda|\leq\epsilon^{-1}}|\hat{\psi}(\lambda)|\geq \eta$. 
Hence for all $j$ such that $2^j\epsilon\geq\epsilon^{-1}$, using~(\ref{eq:fjApprox}), 
$$
\inf_{|\lambda|\leq\epsilon} 2^{-j/2} |H_j(\lambda)| \geq 
\inf_{|\lambda|\leq\epsilon}|\hat{\phi}(\lambda)\hat{\psi}(2^j\lambda)| - C_1 2^{-j\alpha} \sup_{|\lambda|\leq\epsilon} |\lambda|^M
\geq \eta/2 - C_1 2^{-j\alpha} \epsilon^M \; ,
$$
which is positive for $j$ large enough. 
\end{proof}

\section{Some useful inequalities}\label{sec:useful-inequality}

\begin{lemma}\label{lem:empVarL2bound}
Let $p$ be a positive integer. For all $\Cset^p$-valued function $\bg\in L^2(-\pi,\pi)$ and $n\geq1$, define 
\begin{equation}\label{eq:Nn}
M_n(\bg)\eqdef \left\{\sum_{k\in\Zset} \left(1-\frac{|k|}{n}\right)_+ |\bc_k|^2\right\}^{1/2} \;,
\end{equation}
where $|\cdot|$ denotes the Euclidean norm in any dimension and $\bc_k=\int_{-\pi}^\pi \bg(\lambda)\,\rme^{\rmi k
  \lambda}\,d\lambda$. Then, for all $\bg_1$ and $\bg_2$ in $L^2(-\pi,\pi)$,  
\begin{equation}\label{eq:NnLipschitz}
|M_n(\bg_1)-M_n(\bg_2)| \leq \sqrt{2\pi} \, \left(\int_{-\pi}^\pi |\bg_1(\lambda)-\bg_2(\lambda)|^2\,d\lambda\right)^{1/2} \; .
\end{equation}
Moreover, for all $\bg$ in $L^2(-\pi,\pi)$, as $n\to\infty$, 
\begin{equation}\label{eq:NnLim}
M_n(\bg) \to \sqrt{2\pi}\,\left(\int_{-\pi}^\pi |\bg(\lambda)|^2\,d\lambda\right)^{1/2} \; .
\end{equation}
\end{lemma}
\begin{proof}
Suppose $p=1$ (the proof for $p>1$ is identical).
Observe that the RHS in~(\ref{eq:Nn}) is a norm on $(c_k)_{k\in\Zset}\in l^2(\Zset)$ which is bounded by the $l^2$ norm
$\left(\sum_{k\in\Zset}|c_k|^2\right)^{1/2}$. Thus, by Parseval Theorem, $M_n(g)$ is a norm on $g\in L^2(-\pi,\pi)$ which is bounded by
$\sqrt{2\pi}\|g\|_2$. Hence $|M_n(g_1)-M_n(g_2)| \leq M_n(g_1-g_2) \leq \sqrt{2\pi} \, \|g_1-g_2\|_2$. Finally,~(\ref{eq:NnLim}) is
obtained by dominated convergence.
\end{proof}

Denote by $\trace(A)$ and $\rho(A)$ the trace and the spectral radius of a matrix $A$. Recall that $\rho(A)$ is the maximum
of the modulus of the eigenvalues of $A$. 

\begin{lemma}\label{lem:CovRadBound}
Let $\{\xi_\ell,\,\ell\in\Zset\}$ be a stationary process with spectral density $g$
and let $\Gamma_n$ be the covariance matrix of $[\xi_1, \dots, \xi_n]$. Then,
$ \rho(\Gamma_n)\leq 2\pi\,\|g\|_\infty$.
\end{lemma}
\begin{proof}
 Since $\Gamma_n$ is a non-negative definite matrix, $ \rho(\Gamma_n)=\sup_{\bx\in\Rset^n,|\bx|\leq1} \bx^T\Gamma_n\bx$,
 where $|\bx|$ is the Euclidean norm of $\bx$. For all $\bx\in\Rset^n$, we may write
$$
\bx^T\Gamma_n\bx  = 
\int_{-\pi}^\pi g(\lambda)\left|\sum_{\ell=1}^n\bx_\ell\,\rme^{-i\ell\lambda}\right|^2\,d\lambda
\leq \|g\|_\infty \int_{-\pi}^\pi \left|\sum_{\ell=1}^n\bx_\ell\,\rme^{-i\ell\lambda}\right|^2\,d\lambda=2\pi\,\|g\|_\infty\,|\bx|^2 \; .
$$
\end{proof}

\begin{proposition}\label{prop:meanlogquadgauss}
Let $\xi$ is a zero-mean $n \times 1$ Gaussian vector with covariance $\Gamma$.
Then there exists a universal constant $C$ independent of $n$ such that for any $n \times n$ non-negative symmetric matrices
$A$ satisfying $\trace(A \Gamma) > 0$, 
\begin{align}
\label{eq:boundmean}
\left| \PE \left( \log[\xi^T A \xi] \right)  -  \log \left( \PE\left[\xi^T A \xi\right] \right) \right|  & \leq
C \left( 1 \wedge \frac{\rho^{2}(A) \rho^{2}(\Gamma)}{\PVar(\xi^T A \xi) } \right) \; ;\\
\label{eq:boundvar}
\PVar \left( \log[\xi^T A \xi] \right) & \leq C \; .
\end{align}

Let $[\xi^T, \tilde{\xi}^T]^T$ be a zero-mean $(n+\tilde{n})\times1$ Gaussian vector such that $\PCov(\xi)= \Gamma$ and
$\PCov(\tilde{\xi})= \tilde{\Gamma}$. 
Then there exists a universal constant $C$ independent of $n$ and $\tilde{n}$ such that for any $n \times n$ and  $\tilde{n}
\times \tilde{n}$ non-negative symmetric matrices $A$ and $\tilde{A}$ satisfying $\trace(A \Gamma) > 0$ and $\trace(\tilde{A}
\tilde{\Gamma}) > 0$, 
\begin{multline}
\label{eq:boundcov}
\left| \PCov\left( \log[\xi^T A \xi] , \log[\tilde{\xi}^T \tilde{A} \tilde{\xi}] \right)  -
\frac{\PCov( \xi^T A \xi, \tilde{\xi}^T \tilde{A} \tilde{\xi})}
{\PE\left[\xi^T A \xi\right]\PE\left[\tilde{\xi}^T \tilde{A}\tilde{\xi}\right]} 
\right| \leq \\
C
\left\{ \frac{\rho^3(A) \rho^3(\Gamma)}{\PVar^{3/2}(\xi^T A \xi)}
\vee \frac{\rho^3(\tilde{A}) \rho^3(\tilde{\Gamma})}{\PVar^{3/2}(\tilde{\xi}^T \tilde{A} \tilde{\xi})}
\right\}
\eqsp.
\end{multline}
\end{proposition}
\begin{proof}
Let $k$  be the rank of $\Gamma$ and $Q$ be $n\times k$ full rank matrix such that
$QQ^T=\Gamma$. Let $\zeta\sim\calN(0,I_k)$, where $I_k$ is the identity matrix of size $k\times k$. For any unitary matrix
$U$, $U\zeta\sim\calN(0,I_{k})$ and hence $QU\zeta$ has same distribution as $\xi$.
Moreover, since $A$ is symmetric, so is $Q^TAQ$.
We may choose an unitary matrix $U$ such that $ \Lambda\eqdef U^T(Q^TAQ)U$ is a diagonal matrix with non-negative entries. Furthermore,
\begin{equation}
\label{eq:EqDistIndep}
\zeta^T\Lambda\zeta= (QU\zeta)^T A (QU\zeta) \eqdist \xi^T A \xi \; ,
\end{equation}
where $\eqdist$ denotes the equality of distributions.
Since $\Lambda$ is diagonal with non-negative diagonal entries  $(\lambda_{i})_{i=1,\dots,k}$,
$\zeta^T\Lambda\zeta$ is a sum of independent r.v.'s of the form $\sum_{i=1}^{k} \lambda_{i} \zeta_{i}^2$.
Since $\PE\zeta_i^2=1$ and $\PVar(\zeta_i^2)=2$,
we get from~(\ref{eq:EqDistIndep}) that $ \sum_{i=1}^{k}\lambda_{i}=\PE\left[\zeta^T\Lambda\zeta\right]=\PE\left[\xi^T A \xi\right] = \trace(A \Gamma)>0$
and $\PVar \left[\xi^T A \xi\right]=\PVar\left[\zeta^T\Lambda\zeta\right]=2\sum_{i=1}^{k} \lambda_{i}^2$.
Now set
\begin{equation}
\label{eq:Sn}
S \eqdef \frac{\xi^T A \xi}{\PE[\xi^T A \xi]} \eqdist
\sum_{i=1}^{k} d_i \zeta_i^2 \quad \text{with} \quad d_i \eqdef \frac{\lambda_i}{\sum_{j=1}^{k} \lambda_j} \eqsp,
\end{equation}
so that 
\begin{equation}
\label{eq:ESVarS}
\PE[S] = 1 \quad\text{and}\quad \PVar(S) = 2\, \|d\|^2 \; ,
\end{equation}
where $\|d\|^2 \eqdef \sum_{i=1}^{k} d_{i}^2$. The quantities of interest in~(\ref{eq:boundmean}) and~(\ref{eq:boundvar})
become 
\begin{equation}
\label{eq:SnEsp}
\PE \left( \log \left[ \xi^T A \xi \right] \right) - \log \left[ \PE \left( \xi^T A \xi \right) \right] =
\PE \left( \log[S ] \right) 
\end{equation}
and
\begin{equation}
\label{eq:SnVar}
 \PVar(\log \left[ \xi^T A \xi \right])= \PVar \left( \log[S ] \right) \; .
\end{equation}
Since $\sum_{i=1}^{k} \lambda_{i} \geq  \sum_{i=1}^{k} \lambda_{i}^2 /\max_{1 \leq i \leq k} \lambda_{i}$
and $\rho(\Lambda)=\max_{1 \leq i \leq k} \lambda_{i}$, we get
\begin{equation}\label{eq:DntoZero}
\|d\|^2 = \frac{\sum_{i=1}^{k} \lambda_{i}^2}{ \left( \sum_{i=1}^{k} \lambda_{i} \right)^2} \leq
\frac{\rho^2(\Lambda)}{\sum_{i=1}^{k} \lambda_{i}^2}
= \frac{2\rho^2(\Lambda)}{\PVar(\xi^T A \xi)}
\leq \frac{2\rho^2(A)\,\rho^2(\Gamma)}{\PVar(\xi^T A \xi)} \eqsp .
\end{equation}
This is the quantity which appears in~(\ref{eq:boundmean}), and we will therefore express bounds in terms of $\|d\|$. 

Denote by $F$ the distribution function of $S$, that is $F(x)=\prob(S\leq x)$.
Observe that $F(0)=0$ since $S$ is a non-negative weighted sum of independent central chi-squares and that all the weights do not
vanish. To obtain exponential bounds on $F$, observe that, by standard computations on the chi-square distribution, one has, for
$t>-(2 \max_{1 \leq i \leq k} d_i)^{-1}$,   
\begin{equation}\label{eq:carf}
\PE \left[ e^{-t S} \right]=\prod_{i=1} ^{k}\PE \left[ e^{-t d_{i} \zeta_i^2} \right]=
\prod_{i=1} ^{k} (1+2 d_i t)^{-1/2} \; .
\end{equation}
Therefore, for any $t>0$ and $x>0$,
\begin{align}\label{eq:ExpoBound0}
\log \left[ F(x) \right] \leq \log \left[ \rme^{xt} \PE(\rme^{-t S}) \right]  = xt-(1/2) \sum_{i=1} ^{k} \log(1+2 d_{i} t)
\eqsp .
\end{align}
Using~(\ref{eq:ExpoBound0}), we derive two bounds for $F(x)$ by choosing $t$ adequately. One bound, which will not depend on
$\|d\|$ is for $x$ around 0, the other one, which will improve as $\|d\|$ decreases, is for $x$ in $(0,1/2)$.  

To get the first bound, observe that, for $t \geq 0$, $ \prod_{i=1} ^{k} (1 + 2 d_{i}  t) \geq 1 + 2 t \sum_{i=1} ^{k} d_{i}
= 1+ 2t $, $ \sum_{i=1} ^{k} \log(1+2 d_{i} t) \geq \log(1+2t)$. Plugging this inequality in~(\ref{eq:ExpoBound0}) and setting
$t=1/(2x)$ yields 
\begin{equation}
\label{eq:expobound2}
F(x)\leq e^{1/2} \left( \frac{x}{1+x}  \right)^{1/2}\leq \sqrt{\rme \, x},\quad x>0  \eqsp.
\end{equation}
Let $p\geq1$ and $\alpha\in\Rset$. Since $\lim_{x \to 0^+} |\log(x)|^pF(x)= 0$, 
integration by parts and \eqref{eq:expobound2} give that, 
$$
\int_0^1|\log(x)|^p dF(x) = 
p\int_0^1 |\log(x)|^{p-1} x^{-1} F(x) dx \leq \rme^{1/2}\int_0^1|\log(x)|^{p-1}x^{-1/2} dx \; ,
$$
which is a finite constant. Since $\sup_{x\geq1}(|\log(x)|^p/x)$ is finite and $\PE S=1$, we get that
$\int_1^\infty |\log(x)|^p dF(x)$ is bounded by a constant only depending on $p$ and thus  
\begin{equation}
\label{eq:UnifBound1}
\PE[|\log S|^p]=\int_0^1|\log(x)|^p dF(x) + \int_1^\infty |\log(x)|^p dF(x) \leq K_p\; ,
\end{equation}
where $K_p$ is constant only depending on $p$. This bound proves the left part of the $\wedge$ sign in (\ref{eq:boundmean}). 

We now derive a second bound on $F(x)$ which will yield the right part of the $\wedge$ sign in (\ref{eq:boundmean}).
Since the second derivative of $\log(1+u)$ has absolute value at most 1 for all $u\geq0$, we have, by Taylor's formula,
that, for any $t \geq 0$, $\log \left( 1+2 d_{i} t \right) \geq 2 d_{i} t -2 d^2_{i} t^2$. Applying this
to~(\ref{eq:ExpoBound0}) and using $\sum_{i=1}^k d_i=1$, we get
$$
\log [ F(x) ] \leq \left( x- 1 \right)t  +  t^2 \|d\|^{2}, \quad x>0 \; .
$$
Setting this time $t= \|d\|^{-1}$, we obtain the following exponential bound:
\begin{equation}
\label{eq:expobound1}
F(x)\leq \exp \left[ - (1-x) \|d\|^{-1} +1 \right]\leq \exp \left[ -\|d\|^{-1}/2 +1 \right],\quad x\in(0,1/2) \eqsp.
\end{equation}
Using the relation $a\wedge b\leq\sqrt{ab}$, $a,b\geq0$, we can combine~(\ref{eq:expobound2}) and~(\ref{eq:expobound1}) to
get
$$
F(x) \leq x^{1/4}\exp( -\|d\|^{-1}/4 + 3/4),\quad  x\in(0,1/2) \eqsp.
$$
With this last bound of $F$ at hand, we can improve the bound established in~(\ref{eq:UnifBound1}) as follows.
Let $p\geq1$. Since $|\log(x)|^px^{1/4}$ is bounded on $x\in(0,1/2)$ and 
$|\log(x)|^{p-1}x^{-3/4}$ is integrable  on $x\in(0,1/2)$,   we have, by integration by parts,
\begin{align}
\nonumber
\PE \left( |\log(S)|^p \1 \{ S \leq 1/2 \} \right) & 
\leq \left[|\log(x)|^p\,F(x)\right]_0^{1/2}+p\int_0^{1/2}|\log(x)|^{p-1}x^{-1}\,F(x)\,dx \\
\label{eq:bound:lower_tail}
&\leq C_p\,\exp(-\|d\|^{-1}/4) \leq C_{p,\alpha}\, \|d\|^\alpha \; ,
\end{align}
where $C_p$ and $C_{p,\alpha}$ are constants only depending on $p$ and $(p,\alpha)$. 
For $x\in(0,1)$, we have $|\log(x) - (x-1)|\leq |\log(x)|$, and for  $x\geq1/2$, a Taylor expansion gives $|\log(x) - (x-1)|\leq 2(x-1)^2$,
since the second derivative of $\log(x)$ has absolute value at most $(1/2)^{-2}=4$. Hence, for any $x > 0$,
\begin{equation}
\label{eq:borne-elementaire}
| \log(x) - (x-1)| \leq   |\log(x)| \1_{[0,1/2]}(x) + 2 \, (x-1)^2 \1_{[1/2,\infty)}(x) \eqsp .
\end{equation}
Since $\PE[S-1] =0$ and $\PE[ (S-1)^2] = \PVar (S)= \|d\|^2\PVar(\zeta_1^2)= 2 \|d\|^2$,  using~(\ref{eq:borne-elementaire})
and~\eqref{eq:bound:lower_tail} with $p=1$ and $\alpha=2$, we get
\begin{align*}
\left| \PE\left[ \log(S) \right] \right| = \left| \PE\left[ \log(S) - (S-1) \right] \right|
&\leq  \PE \left[ \left| \log(S)  \right| \1 \{ |S | \leq 1/2 \} \right]  + 2 \PE[ (S-1)^2]  \\
&\leq  C_{1,2} \|d\|^2  + 4 \|d\|^2 \eqsp .
\end{align*}
Applying~(\ref{eq:SnEsp}) and~(\ref{eq:DntoZero}), we get the inequality~(\ref{eq:boundmean}) with the right part of the
$\wedge$ sign.  

The bound~(\ref{eq:boundvar}) is obtained by applying~(\ref{eq:UnifBound1}) since, by~(\ref{eq:SnVar}),
$\PVar \left( \log[\xi^T A \xi] \right) = \PVar \left( \log(S) \right) \leq \PE\log^2 S$.

We now prove \eqref{eq:boundcov}. Define $\tilde{k}$, $\tilde{d}_{i}$ and $\tilde{S}$ as we did ${k}$,
${d}_{i}$ and ${S}$.  The LHS in~(\ref{eq:boundcov}) then  reads
\begin{multline*}
\PE \left[ \log(S) \log(\tilde{S}) \right] - \PE \left[ (S-1) (\tilde{S} -1) \right]  =
\PE \left[ (S-1) ( \log(\tilde{S}) - (\tilde{S}-1) ) \right] \\
+ \PE \left[ (\tilde{S}-1) ( \log(S) - (S-1) ) \right] +
\PE \left[ \left( \log(S) - (S-1) \right) ( \log(\tilde{S}) - (\tilde{S} -1)  )\right] \eqsp.
\end{multline*}
We will provide a bound for the first term of the RHS of this display, the other terms being treated similarly.
By using \eqref{eq:borne-elementaire} and the Cauchy-Schwarz inequality,
\begin{align*}
&\left|\PE \left[ (S-1) \left( \log(\tilde{S}) - (\tilde{S}-1) \right) \right] \right| \\
& \quad \leq
\PE \left[ (S-1)  | \log( \tilde{S}) | \1_{[0,1/2]}(\tilde{S}) \right] 
+ 2 \, \PE \left[ (S-1) (\tilde{S}-1)^2 \right] \\
& \quad \leq
\left( \PE |S-1|^2 \,\PE \left[  | \log( \tilde{S}) |^{2} \1_{[0,1/2]}(\tilde{S})  \right] \right)^{1/2}  
+ 2 \, \left( \PE |S-1|^2 \,\PE |\tilde{S}-1|^{4} \right)^{1/2} \; .
\end{align*}
In view of~(\ref{eq:DntoZero}), it remains to show that the two last terms are $O(\|d\|^3 \vee \|\tilde{d}\|^3)$.
By definition, $\tilde{S}-1= \sum_{i=1} ^{\tilde{k}} \tilde{d}_{i}  (\zeta_{i}^2-1)$, where $\{ \zeta_{i}  \}_{1 \leq k \leq \tilde{k}}$
are i.i.d. standard normal. Therefore,
\[
 \PE |\tilde{S}-1|^{4}  = \sum_{i=1} ^{\tilde{k}} \tilde{d}_{i}^4 \mathrm{cum}_{4}(\zeta_{1}^2)
+ 3 \left( \sum_{i=1} ^{k} \tilde{d}_{i}^2 \right)^2 \PVar^2(\zeta_{1}^2),
\]
where $\mathrm{cum}_4(Z)$ is the fourth-order cumulant of the random variable $Z$. Since $\sum_{i=1} ^{\tilde{k}} \tilde{d}_{i}^4 \leq
\left ( \sum_{i=1}^{\tilde{k}} \tilde{d}_{i}^2 \right)^2$, we obtain  that $\PE |\tilde{S}-1|^{4} \leq C \|\tilde{d}\|^4$
for some constant $C$. Therefore,
\[
\left( \PE |S-1|^2  \right)^{1/2} \left( \PE |\tilde{S}-1|^{4} \right)^{1/2}  
\leq C \|d\| \|\tilde{d}\|^2 \leq C (\|d\|^3 \vee \|\tilde{d}\|^3) \eqsp.
\]
Applying~\eqref{eq:bound:lower_tail} with $p=2$ and $\alpha=4$, we have
\[
\left( \PE |S-1|^2 \right)^{1/2} \left( \PE \left[ | \log( \tilde{S}) |^{2} \1_{[0,1/2]}(\tilde{S})  \right] \right)^{1/2}
\leq(\sqrt{2}\|d\|)\,(C_{2,4}\,\|\tilde{d}\|^2)
\leq \sqrt{2}\,C_{2,4} \, (\|d\|^3 \vee \|\tilde{d}\|^3) \eqsp,
\]
for some constant $C'$, which concludes the proof.
\end{proof}


\small

\bibliography{\BIBDIR/lrd}
\bibliographystyle{ims}
\end{document}